\documentclass[12pt]{amsart}
\usepackage{amsmath, amssymb}
\usepackage{mathrsfs}
\newcommand{\no}[1]{#1}
\renewcommand{\no}[1]{}
\no{\usepackage{times}\usepackage[subscriptcorrection, slantedGreek, nofontinfo]{mtpro}
\renewcommand{\Delta}{\upDelta}}
 \usepackage{color}
\usepackage{enumerate,comment}
\usepackage{amssymb}
\usepackage{graphicx}
\usepackage{dsfont}
\usepackage{epstopdf}
\graphicspath{{./pics/}}
\usepackage{booktabs}
\usepackage{threeparttable}
\usepackage{marginnote}
\usepackage{xcolor}
\usepackage{marginnote}
\usepackage{bm}
\usepackage{algorithm}
\usepackage{algpseudocode}
\usepackage{threeparttable}
\usepackage[hidelinks]{hyperref}

\date{\today}
\newcommand{\bel}{\begin{equation} \label}
\newcommand{\ee}{\end{equation}}
\def\beq{\begin{equation}}
\def\eeq{\end{equation}}
\newcommand{\bea}{\begin{eqnarray}}
\newcommand{\eea}{\end{eqnarray}}
\newcommand{\beas}{\begin{eqnarray*}}
\newcommand{\eeas}{\end{eqnarray*}}

\newcommand{\R}{\mathbb{R}}

\newtheorem{theorem}{Theorem}[section]

\newtheorem{example}{Example}[section]
\numberwithin{equation}{section}

\renewcommand{\d}{\mathrm{d}}
\allowdisplaybreaks

\def\phi {\varphi}

 \title[]
{Pollutant Concentration Reconstruction in Porous Media from Localized Measurements}
\author{}
\date{}
\begin{document}
\author[A. Blouza, S. Cen and Y. Kian]{Adel Blouza$^1$, Siyu Cen$^1$ and Yavar Kian$^1$}
\begin{abstract}
Reliable environmental monitoring and groundwater protection heavily rely on the ability to reconstruct full pollutant concentration fields from sparse, partial field measurements. In porous media, dissolved contaminants are transported by groundwater flow while undergoing hydrodynamic dispersion and complex chemical or biological reactions. To address practical contamination challenges, this work investigates the reconstruction of pollutant distributions using water samples collected at discrete locations in unsaturated or saturated soils via sampling tubes equipped with ceramic tips. While these suction samplers provide critical local data, they inherently capture only a fraction of the subsurface environment, leaving wide areas unmonitored. Motivated by this real-world monitoring procedure, we propose steady and time-dependent advection-diffusion-reaction models to simulate contaminant transport. These models account for localized pollution sources modeled either as multiple point sources  or as more general path-supported measure-valued sources, which frequently represent leaking infrastructure or linear contamination pathways. Our main contribution is a generalized unique continuation result demonstrating that, under appropriate physical and mathematical assumptions, two concentration fields coinciding on a non-empty observation region must match throughout the entire space-time domain, even when originating from different source configurations. We investigate this problem from both theoretical and numerical perspectives, backing our mathematical uniqueness results with numerical simulations that validate the efficiency and practical utility of our approach. The theoretical proofs combine regularity properties away from source singularities with unique-continuation principles for parabolic and elliptic partial differential equations, providing a robust framework for environmental risk assessment and site remediation.
\end{abstract}
\maketitle

\renewcommand{\thefootnote}{\fnsymbol{footnote}}
\footnotetext{\hspace*{-5mm} 
\begin{tabular}{@{}r@{}p{16cm}@{}}
$^*$
& This work was partially supported by the multidisciplinary project  FORECAST of the university of Rouen Normandy involving Mathematics and Geo-science.
The work of Y. Kian and S. Cen is supported by the French National Research Agency ANR and Hong Kong RGC Joint Research Scheme for the project IdiAnoDiff (grant ANR-24-CE40-7039). \\
$^1$ 
& Univ Rouen Normandie, CNRS, Normandie Univ, LMRS UMR 6085, F-76000 Rouen, France  ( \texttt{adel.blouza@univ-rouen.fr}, \texttt{yavar.kian@univ-rouen.fr}, \texttt{siyu.cen@univ-rouen.fr}).
\end{tabular}}
\section{Introduction}

The identification and localization of pollutant sources in porous media represent critical challenges in environmental modeling, groundwater protection, and pollution management. These issues are particularly acute in coastal aquifers, where climate change and mounting anthropogenic pressures can severely degrade water quality. Once introduced into the subsurface, dissolved pollutants are transported by groundwater flow and dispersed via hydrodynamic processes. Furthermore, chemical and biological reactions can significantly modify their concentrations through degradation and attenuation mechanisms~\cite{sch}. Reliable estimation of pollutant distributions and the precise identification of contamination sources are therefore essential for comprehensive environmental risk assessment, continuous groundwater monitoring, and the design of effective remediation strategies.

The mathematical modeling of groundwater flow and contaminant transport in porous media possesses a rich history. In his seminal monograph, Bear~\cite{Bear1972} established a robust framework describing fluid flow through porous materials alongside the fundamental mechanisms governing the transport and dispersion of dissolved substances (see also~\cite{huyak} for related literature). This framework remains the cornerstone for modeling subsurface contaminant migration. Specifically, advection--diffusion--reaction equations are widely employed to capture the coupled effects of advective groundwater motion, hydrodynamic dispersion, and attenuation reactions.

The inverse characterization of aquifers using indirect measurements has likewise attracted substantial scholarly attention. For instance, Aliouache et al.~\cite{Aliouache2025} (see also~\cite{theis,Yeh2000,Zhu2005} for earlier analyses) investigated the reconstruction of aquifer hydraulic properties, notably transmissivity and storativity, via hydraulic tomography utilizing oscillatory, transient, and steady-state pumping tests. Their methodology leverages the Principal Component Geostatistical Approach to estimate spatial property distributions from hydraulic-head observations, highlighting the critical role of measurement design and aquifer heterogeneity. 

Although the present study shares the overarching objective of recovering unobserved parameters from incomplete data, it addresses a distinct inverse problem. Rather than inferring hydraulic properties from head measurements, we investigate the mathematical determination of pollutant concentration fields and the localization of pollution sources directly from concentration observations.

\begin{figure}[!ht]
  \centering
  \includegraphics[width=0.7\textwidth]{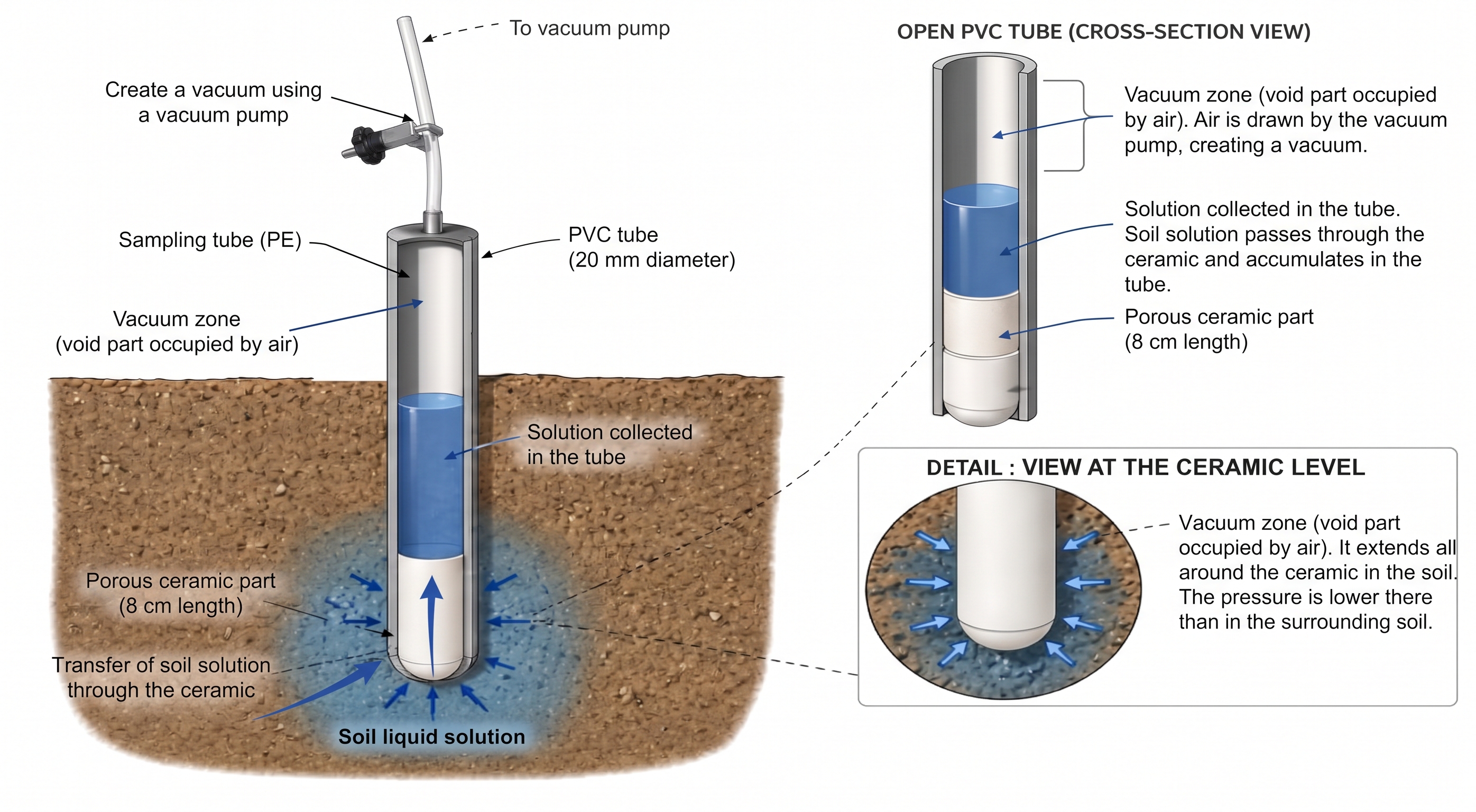}
  \caption{Sampling tube equipped with a ceramic tip for extracting
  water from wet soil.}
  \label{fig:sampling-tube}
\end{figure}

\begin{figure}[!ht]
     \centering	
     \begin{tabular}{cc}
		\includegraphics[width=0.25\textwidth]{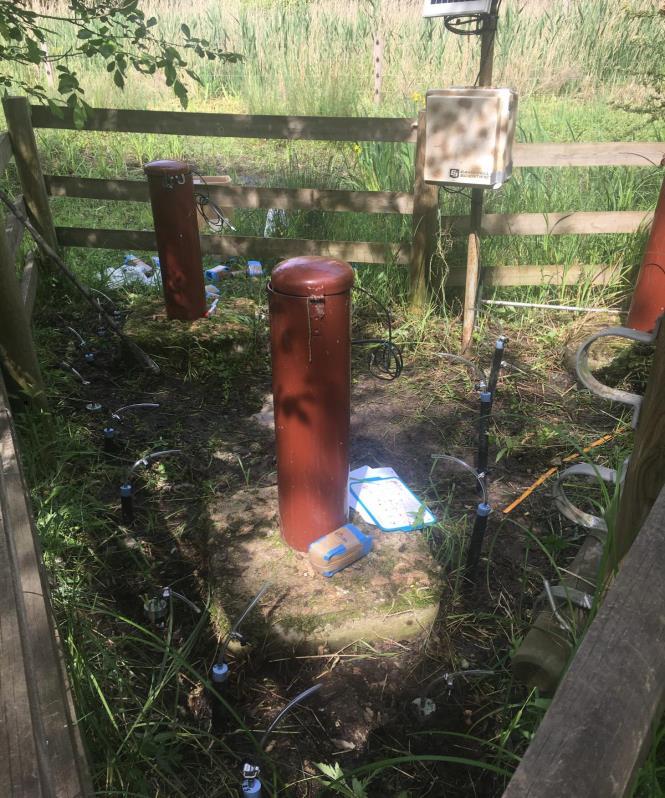}&
		\includegraphics[width=0.5\textwidth]{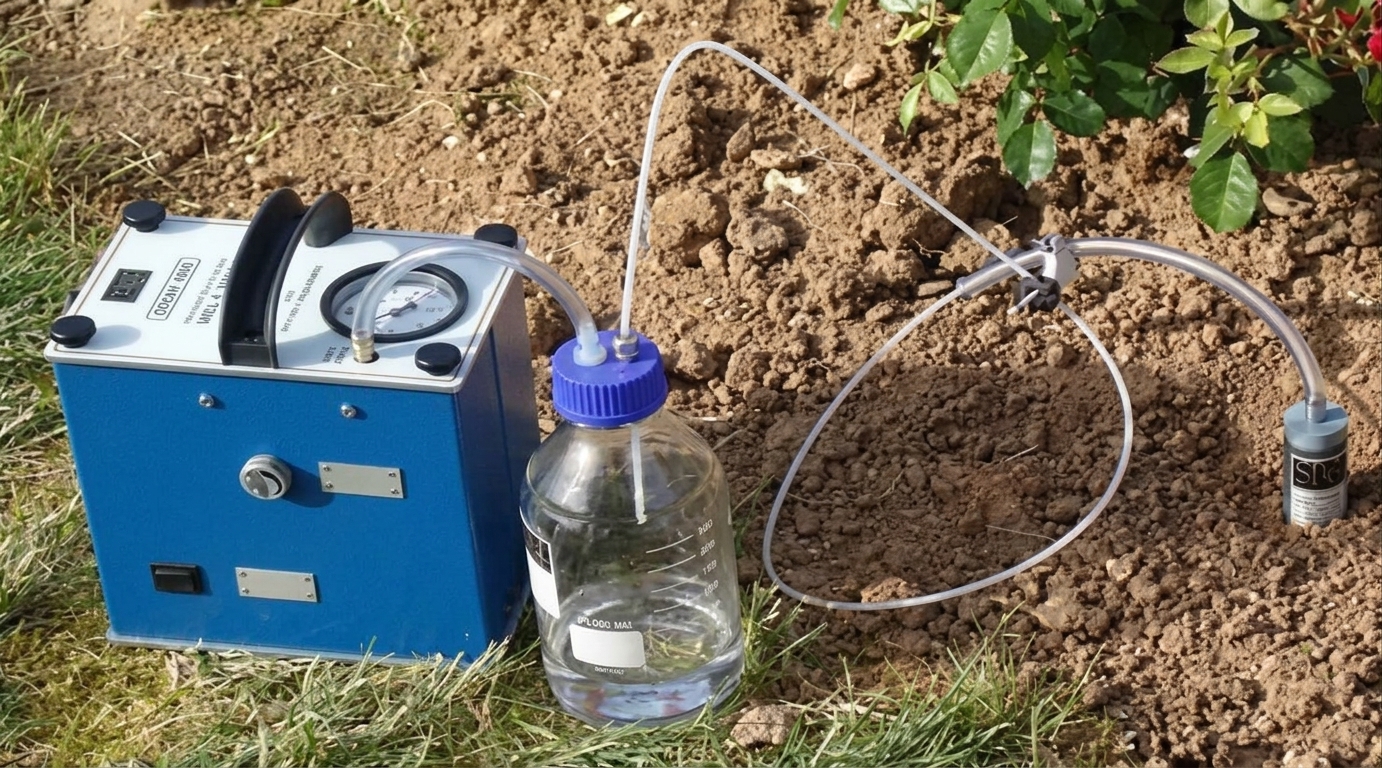} \\
        (a) & (b)
	\end{tabular}  
    \caption{Field sampling setup. (a) Experimental site with multiple sampling tube. (b) Water sampling equipment used for collecting pore water
  from the experimental site. }	\label{fig:exp}
 \end{figure}

The experimental motivation for this work is a monitoring procedure in which water is extracted from wet soil using sampling tubes fitted with ceramic tips, see  Figures~\ref{fig:sampling-tube} and \ref{fig:exp} below (see also \href{https://www.hydroscan.cnrs.fr}{https://www.hydroscan.cnrs.fr} for more details). The collected water samples are subsequently analyzed to determine pollutant concentrations at selected locations within the porous medium. This procedure provides localized information about contamination without requiring direct access to the entire domain. The resulting measurements may be spatially sparse and may not cover all regions affected by pollutant transport. Consequently, a fundamental mathematical question arises: can sparse information obtained from a restricted observation region uniquely determine the concentration field across the broader domain?

Let $\Omega\subset\mathbb{R}^{n}$, $n=2,3$, be a bounded domain representing the porous medium, and let $T>0$ denote the final observation time. We consider the concentration $\mathcal{C}(x,t)$ of a dissolved pollutant, at position $x\in\Omega$ and time $t\in(0,T)$, governed by the advection--diffusion--reaction equation
\begin{equation}
\label{AD}
\partial_t\mathcal{C}
-\nabla\cdot(D\nabla\mathcal{C})
+V\cdot\nabla\mathcal{C}
+q\mathcal{C}
=F(x,t),
\qquad (x,t)\in\Omega\times(0,T),
\end{equation}
where $D$ denotes the hydrodynamic diffusion or dispersion coefficient, $V$ is the groundwater velocity field, $q$ is a reaction or degradation coefficient, and $F$ represents the source of the pollutant modeled by a measure-valued source. 

In many environmental applications, pollutant releases occur at localized sites and can be effectively approximated by measure valued sources \cite{w12092415,MS}. Such sources are naturally represented by finite signed Radon  measures, leading to partial differential equations with singular right-hand sides. In the time-dependent setting, we consider a finite collection of $N\in\mathbb{N}:=\{1,2,\ldots\}$ measure-valued sources and write
\begin{equation}
\label{source}
F(x,t) = \sum_{j=1}^{N} \lambda_j(t)\mu_j(x),
\end{equation}
where, for $j=1,\ldots,N$, $\lambda_j\in L^2(0,T)$ represents the time-dependent intensity of the source, and the finite signed Radon measure $\mu_j$ defines its spatial location. In the present article, we consider two different scenarios: 1) $F$ is described by multiple point sources, with $S_j\in\Omega$ being the location of the $j$-th source and $\mu_j=\delta_{S_j}$ the Dirac measure concentrated at $S_j$; 2) in dimension three, the $j$-th Radon measure-valued source $\mu_j$ is supported by a curve. Each of these scenarios covers different practical situations related to our applications of global recovery of pollutant concentration. Since the source locations and intensities are generally unknown, the inverse problem consists in determining the global concentration field and, where feasible, recovering the underlying source configuration
\[
\mathbf{(\bm{\lambda},\bm{\mu})} = ((\lambda_1,\mu_1), \ldots, (\lambda_N,\mu_N))
\]
from partial concentration observations. Depending on the formulation of the inverse problem, the source intensities may either be prescribed or treated as additional unknowns. The singular nature of these sources, coupled with the limited spatial coverage of the measurements, raises fundamental questions regarding uniqueness and identifiability.

To formalize the observation setting, assume that observations are made by $M\in\mathbb{N}$ sampling tubes with ceramic tips, and let $\omega_1, \ldots, \omega_M \subset \Omega$ be $M$ open not-empty sets modeling the observation subdomains where the sampling tubes extract water (see Figure~\ref{fig:sampling-tube} and Figure~\ref{fig:exp}(a)). These observation sets form a finite collection of observation subdomains, defining the composite observation region as
\begin{equation}
\label{omega}
\omega = \bigcup_{i=1}^{M} \omega_i \subset \Omega.
\end{equation}
 In the experimental context of \cite{Aliouache2025, Yeh2000, Zhu2005}, water samples are collected at discrete locations associated with these regions. Each sample provides local information regarding the pollutant concentration in the extracted fluid, subject to the specific sampling protocol and the hydraulic characteristics of the medium. Consequently, such physical measurements are inherently discrete in space. For our theoretical analysis, however, we adopt an idealized observation framework in which the concentration field is assumed to be known continuously on $\omega \times (0,T)$. This formulation allows us to rigorously investigate whether spatially partial concentration data can uniquely determine the field throughout the entire domain. The transition from discrete physical samples to continuous observations over spatial subdomains represents a distinct reconstruction phase that typically necessitates spatial interpolation, auxiliary measurements, or an appropriate statistical model.

This reconstruction problem is closely related to a generalized unique continuation (GUC) property for elliptic and parabolic PDEs, which can be stated as follows:
\begin{itemize}
\item[{\bf(GUC)}] {\em
Two solutions, potentially driven by unknown sources, of the form \eqref{source}, as well as unknown initial and boundary conditions, that agree on an observation region $\omega$ must necessarily coincide elsewhere in the domain $\Omega$.
}
\end{itemize}

Compared with classical unique continuation properties \cite{GL, SS1, SS} (or related data assimilation type of problems \cite{Bu1,Bu2}), the \textbf{(GUC)} property does not concern solutions to homogeneous equations, but rather a suitable class of solutions to equations with measure-valued sources of the form \eqref{source}. Moreover, when stated for equations of the form \eqref{AD} with measure-valued sources of the form \eqref{source}, the \textbf{(GUC)} shares similarities with inverse source problems, where the objective is to recover the source from internal or boundary observations \cite{Andrle2012moving_pollution_sources, DCK, ElBadiaHaDuong2002pollution_detection_problem, GZZ, gong2026identification}. However, the \textbf{(GUC)} is formulated in a more general setting than the existing literature on inverse source problems, where initial and boundary conditions are typically assumed to be known. Classical inverse source problems in such settings can thus be seen as a byproduct of the \textbf{(GUC)}. Indeed, the \textbf{(GUC)} property plays a pivotal role in inverse problems, data assimilation, and environmental monitoring, where the physical state variable is accessible exclusively through partial observations \cite{GK,w12092415, MS}.


The primary objective of this work is to investigate the reconstruction of pollutant concentration fields from partial observations. We examine the \textbf{(GUC)} principle for both steady-state and time-dependent advection-diffusion-reaction models. In the steady-state setting, the concentration is governed by an elliptic equation featuring localized sources, whereas in the time-dependent setting, we analyze the parabolic equation \eqref{AD} incorporating the source term \eqref{source}. Our central theoretical question is whether knowledge of the concentration over $\omega \times (0,T)$ uniquely determines the field throughout $\Omega \times (0,T)$.

The main contribution of this article is a \textbf{(GUC)} result tailored to concentration fields generated by measure-valued sources. More precisely, under the assumptions specified in the main theorems, we prove that if two solutions of the advection-diffusion-reaction equation coincide on the observation region $\omega \times (0,T)$, they must coincide throughout the entire space-time domain $\Omega \times (0,T)$. The proof successfully combines the regularity properties of solutions away from source singularities with geometric arguments and unique-continuation principles for elliptic and parabolic equations. This establishes a rigorous theoretical foundation for determining global concentration fields from spatially partial observations, even when the underlying source configurations are initially unknown. 

To complement the theoretical analysis, we also investigate the reconstruction problem numerically. We propose a two-step algorithm: first, the pollution source, the unknown boundary data, and the initial data are reconstructed from noisy local measurements on $\omega$ using a regularized Levenberg--Marquardt method; second, the reconstructed data are used to solve the forward problem \eqref{AD} to recover the pollutant concentration field $\mathcal{C}$ over the entire domain. A series of numerical experiments illustrates the theory in both steady-state and time-dependent regimes, confirming that both the source configuration and the concentration field can be recovered accurately.

The remainder of the paper is organized as follows. We first introduce  in Section \ref{s2} the mathematical framework and state the core assumptions governing the transport equation, the source terms, and the observation region and we state the corresponding \textbf{(GUC)} results in both the steady-state and time-dependent regimes.  In Section \ref{s3}, we subsequently prove the \textbf{(GUC)} results and discuss their implications for recovering concentration fields from their knowledge on the obseration region \eqref{omega}. Finally, in Section \ref{sec:num} we evaluate the relevance of this theoretical framework to localized concentration measurements in porous media and outline promising directions for future research concerning source localization and numerical reconstruction.

  \section{Mathematical model, main theoretical results and discussions}\label{s2}
From now on, we assume that $\Omega\subset \mathbb{R}^n$, $n= 2,3$, is a  bounded, connected, open set with $C^2$ boundary. Below, we present the mathematical statement and the main results concerning the \textbf{(GUC)} for equation \eqref{AD} in both the steady-state and time-dependent regimes.
\subsection{Steady-state  problem}
We first consider stationary solutions of  equations \eqref{AD}  which are time independent. To this end, let  $V\in W^{1,\infty}(\Omega)^n$, $q\in L^\infty(\Omega)$ and  $D=(d_{i,j})_{1\leq i,j\leq n}\in C^1(\overline{\Omega};\R^{n\times n})$. We define the following second-order elliptic operator 
$\mathcal{L}$  on the domain $\Omega$ by
\begin{equation}\label{As}
\mathcal{L} v(x) :=-\sum_{i,j=1}^n \partial_{x_i} \left( d_{i,j}(x) \partial_{x_j} v(x) \right)+V(x)\cdot\nabla_xv(x)+q(x)v(x),\quad  x\in\Omega.
\end{equation}
Throughout this section, we assume   that the diffusion coefficient matrix $d:=(d_{i,j})_{1 \leq i,j \leq n}$ is symmetric, namely, 
$$
d_{i,j}(x)=d_{j,i}(x),\quad i,j=1,\ldots,n,\  x\in\overline{\Omega},
$$ 
and uniformly elliptic: there exists a constant $c >0$  such that  
\begin{equation}\label{ells}
 \sum_{i,j=1}^n d_{i,j}(x) \xi_i \xi_j \geq c |\xi|^2,\quad  \xi=(\xi_1,\ldots,\xi_n) \in \mathbb{R}^n,\ x\in\overline{\Omega}.
\end{equation}
Let   $N\in\mathbb N$, and let  $S_1,\ldots,S_N\in\Omega$, be distinct points representing source locations, with associated nonzero amplitude $\lambda_1,\ldots,\lambda_N\in \R\setminus\{0\}$. Following the discussions and definitions of \cite[Chapter 8]{GT} and \cite[Section 6, Chapter 2]{LM1}, we say that  $u\in L^2(\Omega)$  is a solution in the transposition sense of the equation  
\begin{equation}\label{eq1S}
\mathcal{L}u=\sum_{j=1}^N\lambda_j\delta_{S_j},\quad\textrm{in } \Omega 
\end{equation}
if 
\begin{equation}\label{eq1-defS}
\int_\Omega u(x)\mathcal{L}^\star\phi(x)dx=\sum_{j=1}^N \lambda_j\phi(S_j),\quad \phi\in C^\infty_0(\Omega),
\end{equation}
where $\delta_{S_j}$ is the Dirac distribution concentrated at $S_j$ and $\mathcal{L}^{\star}$  denotes the formal adjoint of \(\mathcal{L}\) defined by
\begin{equation}\label{L*}
\mathcal{L}^\star w(x) :=-\sum_{i,j=1}^n \partial_{x_i} \left( d_{i,j}(x) \partial_{x_j} w(x) \right)-\nabla_x\cdot(V(x)w)(x)+q(x)w(x),\quad  x\in\Omega.
\end{equation}   
\par 
Similarly, for $\mu_j$, $j=1,\ldots,N$, some finite signed Radon measures supported on $\Omega$, we say that 
$u\in L^2(\Omega)$  solves the equation
\begin{equation}\label{eq1mS}
\mathcal{L}u=\sum_{k=1}^N\lambda_k\mu_k,\quad\textrm{in } \Omega
\end{equation}
in the transposition sense if 
\begin{equation}\label{eq1m-defS}
\int_\Omega u(x)\mathcal{L}^\star\phi(x)dx=\sum_{k=1}^N \lambda_k\int_\Omega \phi(x)d\mu_k(x),\quad \phi\in C^\infty_0(\Omega).
\end{equation}

We consider here solutions $u \in L^2(\Omega)$ of \eqref{eq1S} and \eqref{eq1mS} modeling the stationary concentration of pollution within the domain $\Omega$. Without loss of generality, we assume that the measurements are performed in a weakly polluted region. More precisely, we suppose that there exists a large contaminated area, represented by the non-empty open set $\mathcal{V}$ such that $\overline{\mathcal{V}} \subset \Omega$, and that one of the following conditions holds:
\begin{equation}
\label{U1}
\{S_1, \ldots, S_M\} \subset \mathcal{V},
\end{equation}
\begin{equation}
\label{U2}
\bigcup_{j=1}^M \operatorname{supp}(\mu_j) \subset \mathcal{V}.
\end{equation}
Then, assuming that $\omega \subset \Omega \setminus \overline{\mathcal{V}}$ is a non-empty open subset of $\mathbb{R}^n$, we consider the \textbf{(GUC)} problem associated with the steady-state equations \eqref{eq1S} and \eqref{eq1mS}.
\begin{itemize}
\item[{\bf(SGUC)}] {\em
 Determine uniquely the concentration of pollution $u\in L^2(\Omega)$, seen as a solution in the transposition sense of \eqref{eq1S} or \eqref{eq1mS},  on $\Omega$ from its knowledge on $\omega$.
}
\end{itemize}

We provide the following answers to the {\bf(SGUC)} problem for equations \eqref{eq1S} and \eqref{eq1mS}. We start with \eqref{eq1S} where our answer takes the following form. 

\begin{theorem}\label{t2}
For $j=1,2$, let $N^j\in\mathbb{N}$, 
$\lambda_1^j,\ldots,\lambda_{N^j}^j\in\mathbb{R}\setminus\{0\}$, and let
$S_1^j,\ldots,S_{N^j}^j$ be  distinct points of $\Omega$.
Let $u_j\in L^2(\Omega)$ be a solution in the transposition sense of \eqref{eq1S}
corresponding to $N=N^j$ and
\[
(S_k,\lambda_k)=(S_k^j,\lambda_k^j),\qquad k=1,\ldots,N^j.
\]
Assume further that
\begin{equation}\label{t2a}
u_j\in\bigcap_{r\in(0,r_j)}H^2(\Omega_r^j),
\end{equation}
where
\[
r_j:=\operatorname{dist}
\bigl(\{S_1^j,\ldots,S_{N^j}^j\},\partial\Omega\bigr)
\]
and
\begin{equation}\label{om}
\Omega_r^j
:=\Omega\setminus\bigcup_{k=1}^{N^j}\bar B(S_k^j,r),
\quad r\in(0,r_j),
\end{equation}
with $
\bar B(S_k^j,r):=\{y\in\mathbb{R}^n:\ |y-S_k^j|\leq r\}
$. If
\begin{equation}\label{t2b}
u_1(x)=u_2(x)\quad\text{a.e. }x\in\omega,
\end{equation}
then
\begin{equation}\label{t2c}
u_1(x)=u_2(x)\quad\text{a.e. }x\in\Omega.
\end{equation}
\end{theorem}

We now extend this result to a more general class of source terms, namely finite signed Radon measures supported on   curves in a three-dimensional domain. More precisely, under suitable regularity assumptions away from these curves, we show that the same unique continuation property holds for equations of the form \eqref{eq1mS}.

 \begin{theorem}\label{t3}Let $n=3$ and, for $j=1,2$, let $N^j\in\mathbb N$, $\lambda_1^j,\ldots,\lambda_{N^j}^j\in \R\setminus\{0\}$ and $\mu_1^j,\ldots,\mu_{N^j}^j$ be $N^j$ distinct finite signed Radon measure on $\Omega$. For $j=1,2$ and $k=1,\ldots,N^j$, we fix $\gamma_k^j\in C^\infty([0,1];\R^3)$, satisfying $\gamma_k^j(1)=\gamma_k^j(0)$, and we assume that 
\bel{t3a} \textrm{supp}(\mu_k^j)\subset X_k^j:=\{\gamma_k^j(s):\ s\in[0,1]\}\subset\Omega,\quad j=1,2,\ k=1,\ldots, N^j.\ee
Consider also $u_j\in L^2(\Omega)$, $j=1,2$, solving \eqref{eq1mS} in the transposition sense with $N=N^j$, $(\mu_k,\lambda_k)=(\mu_k^j,\lambda_k^j)$, $k=1,\ldots N^j$, and assume that
\bel{t3b} u_j\in \bigcap_{r\in (0,r_j)} H^2(\mathcal U_r^j),\ee
with $r_j=\underset{1\leq k\leq N^j}{\min}\textrm{dist}(X_k^j,\partial\Omega)$ and
\begin{equation}\label{t3c}\mathcal U_r^j=\{x\in\Omega:\ \textrm{dist}(x,X_k^j)>r,\ k=1,\ldots,N^j\},\quad r\in(0,r_j).
\end{equation}
 Then,  the condition \eqref{t2b}
implies \eqref{t2c}.
\end{theorem}

\subsection{Unsteady-state problem}
Fix  $T>0$ and consider $V\in L^\infty(0,T;W^{1,\infty}(\Omega))^n$, $q\in L^\infty(\Omega\times(0,T))$, $D=(d_{i,j})_{1\leq i,j\leq n}\in C([0,T];C^1(\overline{\Omega};\R^{n\times n}))$. We  define the elliptic operator $\mathcal{A}(t)$, $t\in[0,T]$,  on the domain $\Omega$ by
\begin{equation}\label{A}
\mathcal{A}(t) v(x) :=-\sum_{i,j=1}^n \partial_{x_i} \left( d_{i,j}(x,t) \partial_{x_j} v(x) \right)+V(x,t)\cdot\nabla_xv(x)+q(x,t)v(x),\quad  (x,t)\in\Omega\times(0,T).
\end{equation}                                                     
We assume here that the diffusion  coefficient matrix $d:=(d_{i,j})_{1 \leq i,j \leq n}$ is symmetric, i.e., $d_{i,j}(x,t)=d_{j,i}(x,t)$, for
any $ (x,t)\in\overline{\Omega}\times[0,T]$, $i,j = 1,\ldots,n,$ and fulfills the following ellipticity condition
\begin{equation}\label{ell}
\exists c>0 :\,  \sum_{i,j=1}^n d_{i,j}(x,t) \xi_i \xi_j \geq c |\xi|^2,\quad  \xi=(\xi_1,\ldots,\xi_n) \in \mathbb{R}^n,\ (x,t)\in\overline{\Omega}\times[0,T].
\end{equation}

We fix $N\in\mathbb N$,  $S_1,\ldots,S_N\in\Omega$, $\lambda_1,\ldots,\lambda_N\in L^2(0,T)$. Following the definition of \cite[Section 2, Chapter 3]{LM1}, we say that  $\mathcal C\in L^2(\Omega\times(0,T))$  solves in transposition sense the equation
\begin{equation}\label{eq1}
\partial_t\mathcal C+\mathcal{A}(t)\mathcal C=\sum_{j=1}^N\lambda_j(t)\delta_{S_j}(x),\quad \textrm{in } \Omega\times(0,T)
\end{equation}
 if 
\begin{equation}\label{eq1-def}
\int_0^T\int_\Omega \mathcal C(x,t)(-\partial_t\phi+\mathcal{A}^\star(t)\phi)(x,t)dxdt=\sum_{j=1}^N\int_0^T\lambda_j(t)\phi(S_j,t)dt,\quad \phi\in C^\infty_0(\Omega\times(0,T)),
\end{equation}
where
\begin{equation}\label{A*}
\mathcal{A}^\star(t) v(x) :=-\sum_{i,j=1}^n \partial_{x_i} \left( d_{i,j}(x,t) \partial_{x_j} v(x) \right)-\nabla_x\cdot(V(x,t)v)(x)+q(x,t)v(x),\quad  (x,t)\in\Omega\times(0,T).
\end{equation}   
Similarly, for $\mu_j$, $j=1,\ldots,N$, some finite signed Radon measures supported on $\Omega$, we say that 
$\mathcal C\in L^2(\Omega\times(0,T))$  solves in the transposition sense the equation
\begin{equation}\label{eq1m}
\partial_t\mathcal C+\mathcal{A}(t)\mathcal C=\sum_{j=1}^N\lambda_j(t)\mu_j(x),\quad  \textrm{in } \Omega\times(0,T)
\end{equation}
 if 
\begin{equation}\label{eq1m-def}
\int_0^T\int_\Omega \mathcal C(x,t)(-\partial_t\phi+\mathcal{A}^\star(t)\phi)(x,t)dxdt=\sum_{j=1}^N\int_0^T\int_\Omega \lambda_j(t)\phi(x,t)d\mu_j(x)dt,\quad \phi\in C^\infty_0(\Omega\times(0,T)).
\end{equation}

We consider here solutions $\mathcal C\in L^2(\Omega\times(0,T))$ of \eqref{eq1} or \eqref{eq1m} modeling the concentration of pollution in the area $\Omega$ along the time interval $(0,T)$. Similarly to the steady-state problem, we assume $\omega \subset \Omega \setminus \overline{\mathcal{V}}$ is a non-empty open subset of $\mathbb{R}^n$, with $\mathcal V$ satifying \eqref{U1} or \eqref{U2}, and we consider the following \textbf{(GUC)} problem for the unsteady-state equations \eqref{eq1} and \eqref{eq1m}.
\begin{itemize}
\item[{\bf(NSGUC)}] {\em
Determine uniquely the concentration of pollution $\mathcal C\in L^2(\Omega\times(0,T))$, seen as a solution of \eqref{eq1} or  \eqref{eq1m} in the transposition sense,  on $\Omega\times(0,T)$ from its knowledge on $\omega\times(0,T)$.
}
\end{itemize}

Our answers to problem {\bf(NSGUC)} take the following form 

\begin{theorem}\label{t1} For $j=1,2$, let $N^j\in\mathbb N$, $\lambda_1^j,\ldots,\lambda_{N^j}^j\in L^2(0,T)$ be $N^j$ non-uniformely vanishing functions, $S_1^j,\ldots,S_{N^j}^j$ be $N^j$ distinct points of $\Omega$. Consider also $\mathcal C_j\in L^2(\Omega\times(0,T))$, $j=1,2$, solving \eqref{eq1} in the transposition sense with $N=N^j$, $(S_k,\lambda_k)=(S_k^j,\lambda_k^j)$, $k=1,\ldots N^j$, and assume that
\bel{t1a} \mathcal C_j\in \bigcap_{r\in (0,r_j)} H^1(0,T;L^2(\Omega_r^j))\cap L^2(0,T;H^2(\Omega_r^j)).\ee
Then,  the condition
\bel{t1b} \mathcal C_1(x,t)=\mathcal C_2(x,t),\quad \textrm{a.e. }(x,t)\in\omega\times (0,T)\ee
implies that
\bel{t1c}\mathcal C_1(x,t)=\mathcal C_2(x,t),\quad \textrm{a.e. }(x,t)\in\Omega\times (0,T).\ee

\end{theorem}

\begin{theorem}\label{t4} Let $n=3$ , for $j=1,2$,  $N^j\in\mathbb N$  and let the conditions of Theorem \ref{t3} be fulfilled. For $j=1,2$ and $k=1,\ldots,N^j$, let   $\lambda^j_1,\ldots,\lambda_{N^j}^j\in L^2(0,T)$ be  non-uniformly vanishing functions, and $\mu^j_1,\ldots,\mu^j_{N^j}$ be  finite signed Radon measure on $\Omega$ satisfying the conditions stated in Theorem \ref{t3}. Consider also $\mathcal C_j\in L^2(\Omega\times(0,T))$, $j=1,2$, solving \eqref{eq1m} in the transposition sense, with $N=N^j$ and  $(\mu_k,\lambda_k)=(\mu^j_k,\lambda^j_k)$, $k=1,\ldots,N^j$, and assume that 
\bel{t4a} \mathcal C_j\in \bigcap_{r\in (0,r_j)} H^1(0,T;L^2(\mathcal U_r^j))\cap L^2(0,T;H^2(\mathcal U_r^j)).\ee
 Then,  the condition \eqref{t1b}
implies \eqref{t1c}.

\end{theorem}

\subsection{Comments on our results}

To the best of our knowledge, our main results provide the first \textbf{(GUC)} results for elliptic and parabolic equations with measure-valued sources. While Theorems \ref{t2} and \ref{t1} consider such problems for point sources, Theorems \ref{t3} and \ref{t4} provide, in dimension $n=3$, an extension of these results to a more general class of measure-valued sources supported on a curve. These findings extend the classical unique continuation properties \cite{GL, SS1, SS}, which are often stated for homogeneous equations, by accommodating potential measure-valued unknown sources in the analysis of the problem. Furthermore, they extend the framework of inverse source problems \cite{Andrle2012moving_pollution_sources, DCK, ElBadiaHaDuong2002pollution_detection_problem, GZZ, gong2026identification}, which can be seen as a byproduct of Theorems \ref{t2} and \ref{t1}, by allowing both unknown initial and boundary conditions. Moreover, in addition to the source recovery, the results of Theorem  \ref{t2}-\ref{t4},
can also be applied to the recovery of the unknown initial and boundary conditions. Indeed, according to the regularity conditions \eqref{t2a}, \eqref{t3b}, \eqref{t1a} and \eqref{t4a}, the recovery of the unknown initial and boundary conditions  follows from the identities \eqref{t2c} and \eqref{t1c} by applying classical trace theorems.

In addition to their mathematical interest arising from their generality, Theorems \ref{t2}-\ref{t4} provide a valuable toolbox for applications to the reconstruction of aquifer hydraulic properties, such as transmissivity and storativity~\cite{Aliouache2025,Yeh2000,Zhu2005}. Indeed, as established by these theorems, localized information about groundwater contamination on $\omega$, extracted via ceramic-tipped sampling tubes (see Figures~\ref{fig:sampling-tube} and \ref{fig:exp}(a)), makes it possible to derive a complete picture of the contaminant distribution across the entire contaminated domain $\Omega$. Such applications to the global reconstruction of groundwater contaminant distribution are further illustrated by the numerical results presented in Section \ref{sec:num}, which fundamentally rely on Theorems \ref{t2} through \ref{t4}.

Observe that the \textbf{(SGUC)} and \textbf{(NSGUC)} properties of Theorems  \ref{t2}-\ref{t4} are stated in a general setting and are subject only to the regularity conditions \eqref{t2b}, \eqref{t3b}, \eqref{t1a}, and \eqref{t4a}. Roughly speaking, these conditions ensure that the solutions enjoy proper regularity away from the support of the source. Such requirements align with the regularity properties imposed in classical unique continuation frameworks \cite{GL, SS1, SS} and are fulfilled under appropriate smoothness assumptions on the various parameters defining the solutions of equations \eqref{eq1S}, \eqref{eq1mS}, \eqref{eq1}, and \eqref{eq1m} (e.g.,  the initial and boundary conditions). 
This general formulation provides flexibility, allowing our results to be tailored to various settings in concrete applications. It also generalizes the statements found in classical inverse source problems \cite{Andrle2012moving_pollution_sources, DCK, ElBadiaHaDuong2002pollution_detection_problem, GZZ, gong2026identification}, where solutions in the transposition sense of equations \eqref{eq1S}, \eqref{eq1mS}, \eqref{eq1}, and \eqref{eq1m} with homogeneous boundary and initial conditions automatically satisfy the regularity conditions \eqref{t2b}, \eqref{t3b}, \eqref{t1a}, and \eqref{t4a}.

Beyond extending in Theorems \ref{t2} and \ref{t1} the framework of point-source recovery for inverse source problems  \cite{Andrle2012moving_pollution_sources, ElBadiaHaDuong2002pollution_detection_problem, GZZ, gong2026identification},  the results of Theorems \ref{t3} and \ref{t4} appear to be the first to address the recovery of a broader class of  measure-valued sources such as signed Radon measures supported on a curve, in dimension $n=3$, for elliptic and parabolic equations.  

The proofs of \textbf{(SGUC)} and \textbf{(NSGUC)}, stated in Theorems~\ref{t2}--\ref{t4}, rely on an approach that combines the unique continuation property for elliptic and parabolic equations with the smoothing assumptions \eqref{t2b}, \eqref{t3b}, \eqref{t1a}, and \eqref{t4a} imposed away from the support of the source. The central idea is to exploit the fact that in both scenarios, the support of the source has null Lebesgue measure  in $\mathbb{R}^n$. Therefore, our goal is to show that the pollutant concentration can be recovered in any open subset of $\Omega$,  that does not intersect the  support of the source \eqref{source}, from the composite observation region \eqref{omega}. Ultimately, this leads to the complete recovery of the concentration field in the entire set $\Omega$.

\section{Proof of the theoretical results}\label{s3}

\subsection{Steady-state problem}
In this section, we prove Theorem~\ref{t2} and ~\ref{t3}, showing that the observation
condition \eqref{t2b} implies the equality \eqref{t2c}.

\begin{proof}[Proof of Theorem~\ref{t2}]
Assume that \eqref{t2b} holds. We aim to prove \eqref{t2c}. Since $\Omega$ is connected, there exists
$r_3\in(0,\min(r_1,r_2))$ such that, for every
$r\in(0,r_3)$, the set
\[
\Omega_r:=\Omega_r^1\cap\Omega_r^2=\Omega\setminus\left(\bigcup_{j=1}^2\bigcup_{k=1}^{N^j}\bar B(S_k^j,r)\right)
\]
is connected and satisfies
$\Omega_r\cap\omega\neq\emptyset$.
Fix $r\in(0,r_3)$, $j\in\{1,2\}$ and let $\phi\in C_0^\infty(\Omega_r)$. Since
$u_j|_{\Omega_r}\in H^2(\Omega_r)$, integration by parts
yields
\[
\begin{aligned}
0
&=\sum_{k=1}^{N^j}c_k^j\phi(S_k^j)=\int_\Omega u_j(x)\mathcal L^*\phi(x)\,dx=\int_{\Omega_r}u_j(x)\mathcal L^*\phi(x)\,dx=\int_{\Omega_r}\mathcal L u_j(x)\phi(x)\,dx.
\end{aligned}
\]
Indeed, the support of $\phi$ is compactly contained in
$\Omega_r$, so all boundary terms vanish. Moreover,  for $j=1,2$ and $k=1,\ldots,N^j$, $S_k^j\not\in\Omega_r$ so 
$\phi(S_k^j)=0$.

Since $\phi\in C_0^\infty(\Omega_r)$ is arbitrary chosen, we
deduce that
\begin{equation}\label{t2d}
\mathcal L u_j=0
\quad\text{in }\Omega_r,
\end{equation}
in the transposition sense and, consequently, almost
everywhere in $\Omega_r$.

Now set $u:=u_1-u_2$ and, applying \eqref{t2a}, observe that
$u|_{\Omega_r}\in H^2(\Omega_r)$. Furthermore, combining
\eqref{t2b} with \eqref{t2d}, we obtain
\begin{equation}\label{eq2S}
\begin{cases}
\mathcal L u=0 & \text{in }\Omega_r,\\
u=0 & \text{in }\omega\cap\Omega_r.
\end{cases}
\end{equation}
Since $\Omega_r$ is open and connected, and
$\omega\cap\Omega_r$ contains a nonempty open subset,
the unique continuation property for elliptic equations
(see, e.g., \cite[Theorem 1]{SS1}) yields $u=0$ on $\Omega_r$.
Consequently, we have $u_1(x)=u_2(x)$, a.e. $x\in \Omega_r$. Since this property is true  for every $r\in(0,r_3)$, we conclude
that
\[
u_1(x)=u_2(x),\qquad
\textrm{ a.e. }x\in
\bigcup_{r\in(0,r_3)}\Omega_r
=
\Omega\setminus
\left\{
S_1^1,\ldots,S_{N^1}^1,
S_1^2,\ldots,S_{N^2}^2
\right\}.
\]
This clearly implies \eqref{t2c}.
\end{proof}

\begin{proof}[Proof of Theorem~\ref{t3}]

 We assume that \eqref{t2b} is fulfilled and we will show \eqref{t2c}. Let us fix 
\bel{X}\mathcal X:=\Omega\setminus(\bigcup_{j=1}^2\bigcup_{k=1}^{N^j}X_k^j).\ee
For $j=1,2$ and $k=1,\ldots,N^j$, since $\gamma_k^j\in C^\infty([0,1];\R^3)$ with $\gamma_k^j(0)=\gamma_k^j(1)$, the curve $X_k^j$ is a regular closed curve of dimension $1$ drawn in $\R^3$ contained into $\Omega$. Thus, by applying the transversality theory developed in \cite[Chapter 2, pp. 67–84]{Hir} (see also \cite[Lemma 2.2]{Nag}),  the   curves $X_k^j$, $j=1,2$ and $k=1,\ldots,N^j$,  have a co-dimension of $2$ and, removing a finite number of closed sets of co-dimension $2$  from the connected open set $\Omega$ leaves the set connected. Therefore, the set $\mathcal X$ is an open  connected and therefore path-connected set of $\R^3$.
Let $x\in\mathcal X$ and choose $y\in\omega$.
We can define a path $\kappa\in C([0,1];\R^3)$, with $A:=\kappa([0,1])\subset\mathcal X$, such that $\kappa(0)=y$ and $\kappa(1)=x$. Fix
\bel{del}\delta':=\min_{j=1,2}\min_{k=1,\ldots,N^j}\textrm{dist}(A,X_k^j)>0,\quad \delta:=\min (\delta',\textrm{dist}(A,\partial\Omega))>0.\ee
Choosing $r_3=\min(r_1,r_2,\frac\delta4)$, we can find $B_1,\ldots,B_p$, $p\in\mathbb N$, balls of diameter less than  $\frac\delta4$ and $\mathcal O$ an open set with smooth boundary satisfying the following properties
\bel{bb1}A\subset \mathcal B:=\bigcup_{k=1}^pB_k\subset\mathcal O\subset \mathcal U_{r_3}:=\mathcal U_{r_3}^1\cap\mathcal U_{r_3}^2,\ee
\bel{bb2}y\in B_1,\quad B_k\cap B_{k+1}\neq\emptyset,\quad x\in B_p,\quad k=1,\ldots,p-1.\ee
In light of  \eqref{t3b}, for $j=1,2$, we have $u_j|_{\mathcal O}\in H^2(\mathcal O)$. Then, fixing $\phi\in C^\infty_0(\mathcal O)$ recalling that $\mathcal O\cap X_k^j=\emptyset$, $j=1,2$ and $k=1,\ldots,N^j$,
  integrating by parts and repeating the arguments of Theorem~\ref{t2}, we find 
$$\begin{aligned}0=\sum_{k=1}^{N^j}\int_{X_j^i}  \underbrace{\phi(x)}_{=0}d\mu_k^j(x)=\sum_{k=1}^{N^j}\int_\Omega \phi(x)d\mu_k^j(x)=\int_\Omega u_j(x)\mathcal L^*\phi(x)dx&=\int_{\mathcal O} u_j(x)\mathcal L^*\phi(x)dx\\&=\int_{\mathcal O} \mathcal L u_j(x)\phi(x)dx.\end{aligned}$$
 Recalling that in the last identity $\phi\in C^\infty_0(\mathcal O)$ is arbitrary chosen, we deduce that, for almost every $x\in \mathcal O$, we have 
$\mathcal L u_j(x)=0$.
Choosing $u=u_1-u_2$ and applying \eqref{bb1}, we observe that $u|_{\mathcal B}\in H^2(\mathcal B)$ and, in a similar way to Theorem \ref{t3}, we find
\bel{eqq1} \mathcal Lu(x)=  0,\quad \textrm{a.e. }x\in\mathcal B.\ee
In particular, using \eqref{bb2}, we find 
\begin{equation}\label{eq2SS}
\begin{cases}
 \mathcal Lu=  0, & \mbox{in }B_1,\\
 u= 0, & \mbox{on } \omega\cap B_1\neq\emptyset.
\end{cases}
\end{equation}
Combining \eqref{eq2SS} with unique continuation properties for elliptic equations (see e.g. \cite[Theorem 1]{SS1}), we obtain 
$$u_1(x)-u_2(x)=u(x)=0,\quad \textrm{a.e. }x\in B_1.$$
Repeating the above argumentation and applying  \eqref{bb2}, we can prove by iteration that $u=0$ almost everywhere  on $\mathcal B$. In particular we have $u=0$ almost everywhere  on $B_p$ which is a neighborhood of $x$. Since $x$ is arbitrary chosen on $\mathcal X$, we deduce that $u_1-u_2=u=0$ almost everywhere  in $\mathcal X$. Therefore, we have
$$u_1(x)=u_2(x),\quad \textrm{a.e. }x\in\Omega\setminus(\bigcup_{j=1}^2\bigcup_{k=1}^{N^j}X_k^j).$$
Recalling that, for each $j=1,2$ and $k=1,\ldots,N^j$, $X_k^j$ is a smooth curve associated with the map $\gamma_k^j\in C^\infty([0,1];\R^3)$, it is well-known that the path $X_k^j$ has three-dimensional
Lebesgue measure zero. The same will be true for the set 
$$\bigcup_{j=1}^{2}\bigcup_{k=1}^{N}X_k^j$$
and we deduce that  $u_1=u_2$ almost everywhere
in $\Omega$, which proves \eqref{t2c}.
\end{proof}

 \subsection{Unsteady-state problem}  We now turn to the unsteady-state problem stated in Theorem~\ref{t1} and ~\ref{t4} . As in the
steady-state case, the proof relies on first establishing
that the difference between two candidate solutions
satisfies the homogeneous equation away from the support of the sources. The desired uniqueness result then follows from
the unique continuation property for parabolic equations.

\begin{proof}[Proof of Theorem~\ref{t1}]
Assume that \eqref{t1b} holds. We aim to prove
\eqref{t1c}. For sufficiently small $r>0$, let
$\Omega_r:=\Omega_r^1\cap\Omega_r^2$.
Under the geometric assumptions on the excluded
neighborhoods, we may choose
$r_3\in(0,\min(r_1,r_2))$ such that, for every
$r\in(0,r_3)$, the set $\Omega_r$ is open and connected
and satisfies $\Omega_r\cap\omega\neq\emptyset$.
Fix $r\in(0,r_3)$,  $j\in\{1,2\}$ and let $\phi\in C_0^\infty(\Omega_r\times(0,T))$.
Since the support of $\phi$ does not intersect the source
points $S_k^j$, we have
\[
\sum_{k=1}^{N^j}\int_0^T
\lambda_k^j(t)\phi(S_k^j,t)\,dt=0.
\]
Combining this with  the fact that, in view of \eqref{t1a}, we have
$\mathcal C_j\in H^1(0,T;L^2(\Omega_r))\cap L^2(0,T;H^2(\Omega_r))$ and, integrating by parts, we obtain
\[
\begin{aligned}
0=\sum_{k=1}^{N^j}\int_0^T
  \lambda_k^j(t)\phi(S_k^j,t)\,dt
&=\int_0^T\int_\Omega
  \mathcal C_j(x,t)
  \bigl(-\partial_t\phi+
  \mathcal A^\star(t)\phi\bigr)(x,t)\,dx\,dt\\
&=\int_0^T\int_{\Omega_r}
  \mathcal C_j(x,t)
  \bigl(-\partial_t\phi+
  \mathcal A^\star(t)\phi\bigr)(x,t)\,dx\,dt\\
&=\int_0^T\int_{\Omega_r}
  \bigl(\partial_t\mathcal C_j+
  \mathcal A(t)\mathcal C_j\bigr)(x,t)
  \phi(x,t)\,dx\,dt.
\end{aligned}
\]
The last equality follows from integration by parts in
space and time. All boundary terms vanish because
$\phi$ is compactly supported in
$\Omega_r\times(0,T)$.

Since $\phi$ is arbitrary, we deduce that
\begin{equation}\label{t1d}
\partial_t\mathcal C_j+
\mathcal A(t)\mathcal C_j=0
\quad\text{in }\Omega_r\times(0,T),
\end{equation}
in the transposition sense and hence almost everywhere
in $\Omega_r\times(0,T)$.

Now set $\mathcal C:=\mathcal C_1-\mathcal C_2$. By the assumed regularity of $\mathcal C_1$ and
$\mathcal C_2$, we have
\[
\mathcal C|_{\Omega_r\times(0,T)}
\in
H^1(0,T;L^2(\Omega_r))
\cap
L^2(0,T;H^2(\Omega_r)).
\]
Combining \eqref{t1b} with \eqref{t1d}, we obtain
\begin{equation}\label{eq2}
\begin{cases}
\partial_t\mathcal C+\mathcal A(t)\mathcal C=0,
& \text{in }\Omega_r\times(0,T),\\
\mathcal C=0,
& \text{on }(\omega\cap\Omega_r)\times(0,T).
\end{cases}
\end{equation}

Since $\Omega_r$ is open and connected and
$\omega\cap\Omega_r\neq\emptyset$,
the unique continuation property for parabolic equations
(see, e.g., \cite[Theorem 1.1]{SS}) yields $\mathcal C=0
\ \text{in }\Omega_r\times(0,T)$. Consequently,
\[
\mathcal C_1(x,t)=\mathcal C_2(x,t),
\quad \textrm{a.e. }(x,t)\in\Omega_r\times(0,T).
\]
As this equality holds for every $r\in(0,r_3)$, we
deduce that
\[
\mathcal C_1(x,t)=\mathcal C_2(x,t)
\quad\textrm{a.e. }(x,t)\in
\left(
\Omega\setminus
\left\{
S_1^1,\ldots,S_{N_1}^1,
S_1^2,\ldots,S_{N_2}^2
\right\}
\right)\times(0,T).
\]
It follows 
that $\mathcal C_1=\mathcal C_2$ almost everywhere 
in $\Omega\times(0,T)$. This proves \eqref{t1c}.
\end{proof}
 
\begin{proof}[Proof of Theorem~\ref{t4}]
We assume that \eqref{t1b} holds and prove \eqref{t1c}. Recall from  Theorem~\ref{t3},
that the set $\mathcal X$, defined by \eqref{X} is a path-connected open set of $\R^3$.
Let $x\in\mathcal X$ and choose $y\in\omega$.
We can define a path $\kappa\in C([0,1];\R^3)$, with $A:=\kappa([0,1])\subset\mathcal X$, such that $\kappa(0)=y$ and $\kappa(1)=x$. Similarly to Theorem~\ref{t3}, fixing $\delta$, given by \eqref{del}, and $r_3=\min(r_1,r_2,\frac\delta4)$, we can find $B_1,\ldots,B_p$, $p$ balls of diameter less than  $\frac\delta4$ and $\mathcal O$ an open set with smooth boundary satisfying the conditions \eqref{bb1}-\eqref{bb2}. Let
$\phi\in C_0^\infty(\mathcal O\times(0,T))$.
Since $\mathcal O\subset\mathcal X=\emptyset$, we have 
$$
\sum_{k=1}^{N^j}\int_0^T\int_\Omega \lambda_k(t)\varphi(x,t)\d\mu_k(x)\, dt=0,\quad j=1,2. 
$$

Combining this with  the fact that, in light of \eqref{t4a} and \eqref{bb1}, we have
$$\mathcal C_j|_{\mathcal O\times(0,T)}\in H^1(0,T;L^2(\mathcal O))\cap L^2(0,T;H^2(\mathcal O))$$ and, integrating by parts in a similar way to Theorem~\ref{t1}, we find
\[
\begin{aligned}
0
=
\int_0^T\int_\Omega
\mathcal C_j(x,t)
\bigl(-\partial_t\phi+
\mathcal A^\star(t)\phi\bigr)(x,t)\,dx\,dt
&=
\int_0^T\int_{\mathcal O}
\mathcal C_j(x,t)
\bigl(-\partial_t\phi+
\mathcal A^\star(t)\phi\bigr)(x,t)\,dx\,dt\\
&=
\int_0^T\int_{\mathcal O}
\bigl(\partial_t\mathcal C_j+
\mathcal A(t)\mathcal C_j\bigr)(x,t)
\phi(x,t)\,dx\,dt.
\end{aligned}
\] Since $\phi\in C_0^\infty(\mathcal O\times(0,T))$ is arbitrary chosen, applying \eqref{bb1}, we conclude that
\begin{equation}\label{t4d}
\partial_t\mathcal C_j+
\mathcal A(t)\mathcal C_j=0
\quad\text{in }\mathcal B\times(0,T)\subset\mathcal O\times(0,T),
\end{equation}
in the transposition sense, and hence almost everywhere.
Now set $\mathcal C=\mathcal C_1-\mathcal C_2$. Then
$\mathcal C|_{\mathcal O\times(0,T)}
\in
H^1(0,T;L^2(\mathcal O))
\cap
L^2(0,T;H^2(\mathcal O))$ and, 
combining \eqref{t1b} with \eqref{t1b} and \eqref{t4d}, we obtain
\begin{equation}\label{eq2m}
\begin{cases}
\partial_t\mathcal C+\mathcal A(t)\mathcal C=0,
&\text{in }B_1\times(0,T),\\
\mathcal C=0,
&\text{in }(\omega\cap B_1)\times(0,T).
\end{cases}
\end{equation}
On the other hand, in view of \eqref{bb2}, the set $B_1\cap\omega$ is an open nonempty subset of $\R^3$ and  the unique
continuation property for parabolic equations
(see, e.g., \cite[Theorem 1.1]{SS}) implies that $\mathcal C=0$ almost everywhere  on $B_1\times(0,T)$.  Repeating the above argumentation and applying  \eqref{bb2}, we can prove by iteration that $\mathcal C=0$ almost everywhere  on $\mathcal B\times(0,T)$. In particular we have $\mathcal C=0$ almost everywhere  on $B_p\times(0,T)$ which is a neighborhood of $(x,t)$, $t\in(0,T)$. Since $x$ is arbitrary chosen on $\mathcal X$, we deduce that $\mathcal C_1-\mathcal C_2=\mathcal C=0$ almost everywhere  in $\mathcal X\times(0,T)$. Therefore, we have
$$\mathcal C_1(x;t)=\mathcal C_2(x,t),\quad \textrm{a.e. }(x,t)\in\mathcal X\times (0,T).$$
Then, repeating the argumentation at the end of the proof of Theorem~\ref{t3}, one can check that
$\mathcal C_1=\mathcal C_2$ almost everywhere  in
$\Omega\times(0,T)$, which proves \eqref{t1c}.
\end{proof}
\section{Numerical reconstruction}\label{sec:num}
In this section, we consider the numerical reconstruction of the pollutant concentration field $u$ or $\mathcal{C}$ over the entire domain, as well as the pollution source configuration $(\bm{\lambda},\bm{\mu})$. We consider the pollutant concentration field as solutions of the boundary value problem
\begin{equation}\label{equabv}
\begin{cases}
\mathcal{L}u = \sum_{j=1}^N \lambda_j \mu_j, &\textrm{in } \Omega, \\
u(x) = g(x), & x \in \partial\Omega,
\end{cases}
\end{equation}
and  the initial boundary value problem
\begin{equation}\label{equa}
\begin{cases}
\partial_t \mathcal{C} + \mathcal{A}(t)\mathcal{C} = \sum_{j=1}^N \lambda_j(t)\mu_j(x), &\textrm{in } \Omega \times (0,T), \\
\mathcal{C}(x,t) = g(x,t), & (x,t) \in \partial\Omega \times (0,T), \\
\mathcal{C}(x,0) = u_0(x), & x \in \Omega,
\end{cases}
\end{equation} 
respectively, where the boundary data $g$ and the initial condition $u_0$ are not necessarily known.

The reconstruction is divided into two steps. First, we employ a least-squares-type algorithm to reconstruct the source configuration $(\bm{\lambda},\bm{\mu})$, the unknown boundary condition $g$, and the initial condition $u_0$. Then, we use these reconstructed data to solve the forward problems \eqref{equa} and \eqref{equabv} and recover the concentration field $u$ and $\mathcal{C}$, respectively.

Throughout this work, the forward problems are solved using the Galerkin finite element method (with conforming piecewise linear elements) in space and the backward Euler scheme in time (for the time-dependent problem) \cite{Thomee:2006}. The noisy measurements $z^\delta \in L^2(\omega)$ and $z^\delta \in L^2(\omega \times (0,T))$ are generated by
\begin{equation}\label{eqn:noise}
    z^\delta(x) = u(x)(1 + \delta \xi), \quad x \in \omega, \quad \text{and} \quad z^\delta(x,t) = \mathcal{C}(x,t)(1 + \delta \xi), \quad (x,t) \in \omega \times (0,T),
\end{equation}
for the steady-state problem \eqref{equabv} and the time-dependent problem \eqref{equa}, respectively. Here, $\delta > 0$ denotes the relative noise level, and $\xi$ is distributed uniformly on $[-1,1]$. All experiments are implemented using the \texttt{FreeFem++} software package \cite{Hecht:2012}.

\subsection{Numerical algorithm}\label{subsec:num_alg}
In this part, we present the numerical algorithm for the reconstruction. Throughout this section, we denote by $\bm{\mu}$ the location parameter of the pollution source (whose concrete form is specified in each numerical experiment), by $\bm{\lambda}$ the intensity of the pollution source, by $g$ the boundary condition, and by $u_0$ the initial condition. In the first step, we reconstruct the unknown parameters $(\bm{\mu},\bm{\lambda},g,u_0)$ from the noisy measurement $z^\delta$ in \eqref{eqn:noise}. We employ the Levenberg--Marquardt algorithm \cite{Levenberg:1944,Marquardt:1963}. The description below is given for the time-dependent problem \eqref{equa}; the steady-state case \eqref{equabv} is obtained as a simple special case by dropping the time variable  and the initial condition $u_0$. Specifically, we define a nonlinear operator $F: (\bm{\mu},\bm{\lambda},g,u_0) \mapsto \mathcal{C}(\bm{\mu},\bm{\lambda},g,u_0)|_{\omega\times(0,T)}$, where $\mathcal{C}(\bm{\mu},\bm{\lambda},g,u_0)$ solves the  forward advection--diffusion--reaction problem \eqref{equa} with parameter $(\bm{\mu},\bm{\lambda},g,u_0)$.
 Since $F$ is affine in $\bm{\lambda}$, we   eliminate the intensity $\bm{\lambda}$ from a linear least-squares problem and iterate directly on the remaining parameters. More precisely, starting from an initial guess $(\bm{\mu}^0,g^0,u_0^0)$, for given $\bm{\mu}^k,\,g^k$ and $u_0^k$ we update the intensity analytically from 
\begin{equation}\label{eqn:compute_q}
    \bm{\lambda}^k=\operatorname*{argmin}_{\bm{\lambda}} \frac{1}{2}\|F(\bm{\mu}^k,\bm{\lambda},g^k,u_0^k)-z^\delta \|_{L^2(\omega\times(0,T))}^2.
\end{equation}
Note $F$ is affine in $\bm{\lambda}$, the minimizer is given explicitly by the normal equations. Then we implement the Levenberg--Marquardt iteration only on the remaining parameters: 
\begin{equation}\label{eqn:compute_mu}
    (\bm{\mu}^{k+1} ,g^{k+1},u_0^{k+1})=\arg\min J_k(\bm{\mu} ,g,u_0), 
\end{equation} 
with the quadratic functional 
\begin{align*} 
    J_k(\bm{\mu},g,u_0)=&\frac{1}{2}\|\tilde F(\bm{\mu}^{k},g^{k},u_0^{k})-z^\delta + \partial_{\bm{\mu}} \tilde F(\bm{\mu}^{k},g^{k},u_0^{k})(\bm{\mu}-\bm{\mu}^k) \\
    &\quad+ \partial_{g} \tilde F(\bm{\mu}^{k},g^{k},u_0^{k})(g-g^k) + \partial_{u_0} \tilde F(\bm{\mu}^{k},g^{k},u_0^{k})(u_0-u_0^k)\|_{L^2(\omega\times(0,T))}^2\\
    &+ \frac{\gamma_{\bm{\mu}}^k}{2}|\bm{\mu}-\bm{\mu}^k|^2 +  \frac{\gamma_{g}^k}{2}\|g-g^k\|_{L^2(\partial\Omega\times(0,T))}^2+ \frac{\gamma_{u_0}^k}{2}\|u_0-u_0^k\|_{L^2( \Omega )}^2
\end{align*}
where $\tilde F(\bm{\mu},g,u_0):=F(\bm{\mu},\bm{\lambda}(\bm{\mu},g,u_0),g,u_0)$ is the reduced forward map,  $\gamma_{\bm{\mu}},\, \gamma_{g}$ and $\gamma_{u_0}$ are penalty parameters, chosen according to the specific problem, $|\cdot|$ denotes the Euclidean norm.  The Jacobian $\partial_{\bm{\mu}}F$ is approximated by the componentwise central finite difference scheme
\begin{equation*}
\partial_{\mu_i} \tilde F(\bm{\mu}) \approx (2\epsilon)^{-1}\bigl(\tilde F(\bm{\mu}+\epsilon \bm{e}_i)-\tilde F(\bm{\mu}-\epsilon \bm{e}_i)\bigr),
\end{equation*}
with $\epsilon=1\times10^{-3}$, where $\bm{e}_i$ denotes the $i$-th coordinate vector. The Jacobians $\partial_{g} \tilde F$ and $\partial_{u_0}\tilde F$ are computed exactly, by  solving forward problem for each nodal basis function.  Note $J_k$ is quadratic in $(\bm{\mu},g,u_0)$, the minimizer can be computed explicitly.

Once the parameters $(\bm{\mu},\bm{\lambda},g,u_0)$ are recovered. We perform a forward solver and reconstruct the pollutant concentration fields $\mathcal{C}$. We summarize the numerical reconstruction in Algorithm \ref{alg:num}.

\begin{algorithm}[htbp]
\caption{Reconstruction algorithm for time-dependent problem}\label{alg:num}
\begin{algorithmic}[1]
\Require  Noisy measurement $z^\delta$ and penalty parameters $\gamma_{\bm\mu},\,\gamma_g$ and $\gamma_{u_0}$.
\State  Step 1: Reconstruction of unknown parameters $(\bm{\mu},\bm{\lambda},g,u_0)$.
\State \hspace{1.12cm} Choose initial guess $(\bm{\mu}^0,g^0,u_0^0)$ and set $k=0$.
\State \hspace{1.12cm} Compute $\bm{\lambda}^k$ using \eqref{eqn:compute_q}.
\State \hspace{1.12cm} Compute $\tilde{F}(\bm{\mu}^k,g^k,u_0^k)=F(\bm{\mu}^k,\bm{\lambda}^k,g^k,u_0^k)$.
\State \hspace{1.12cm} Compute the Jacobians $\partial_{\bm{\mu}}\tilde{F}$, $\partial_{g}\tilde{F} $ and $\partial_{u_0}\tilde{F} $.
\State \hspace{1.12cm} Compute the minimum of $J_k$ and update $(\bm{\mu}^{k+1},g^{k+1},u_0^{k+1})$ via \eqref{eqn:compute_mu}.
\State \hspace{1.12cm} Increase $k$ by one and go to line 3.
\State \hspace{1.12cm} Repeat the above procedure until a stopping criterion is satisfied.
\State Step 2: Reconstruction of pollutant concentration fields $\mathcal{C}$.
\State \hspace{1.12cm}  Solve the forward problem   \eqref{equabv} or \eqref{equa}  using  parameters   $(\bm{\mu},\bm{\lambda},g)$ or $(\bm{\mu},\bm{\lambda},g,u_0)$, respectively,  in Step 1.
\end{algorithmic}
\end{algorithm}

\subsection{Numerical experiments}\label{subsec:num_ex}
In this part, we present the numerical experiments using Algorithm \ref{alg:num}. Throughout the section, we take $\Omega=(0,1)^2$ ($\Omega=(0,1)^3$ for the three dimensional case) and $T=1$. We assume the diffusion coefficient is isotropic, i.e., $D=(d_{i,j})_{1\le i,j\le n}=d \,\mathbb{I}$, with $d\in C^1(\overline{\Omega}) $ and $\mathbb{I}$ denoting the $n\times n$ identity matrix. We also assume the reaction term is zero, i.e., $q\equiv0$. The point source is approximated by a Gaussian, i.e. $\delta_{S}(x)\approx\frac{1}{2\pi\sigma^2}e^{-\frac{|x-S|^2}{2\sigma^2}}$, with $\sigma=0.05$, which avoid numerical instability and makes the computed data depend smoothly on $S$. No such regularization is used for the curve source: since $\mu$ is an explicitly parameterized, the source enters the weak formulation as the line integral and is assembled directly. 
The forward problem is solved on a fine grid with  mesh size $h=\frac{1}{100}$ and time step size $\tau=\frac{1}{600}$, and the inverse problem on a coarse grid with $h=\frac{1}{40}$ and $\tau=\frac{1}{300}$. The accuracy of a reconstruction is measured by the relative errors:
\begin{align*}
    &e_{\bm{\mu}} =\frac{|\bm{\mu}^\dagger-\bm{\mu}|}{|\bm{\mu}^\dagger|},     &&e_{\bm{\lambda}}=\frac{|\bm{\lambda}^\dagger-\bm{\lambda}|}{|\bm{\lambda}^\dagger|},      &&&e_{g}=\frac{\|g^\dagger-g\|_{L^2(\partial\Omega\times(0,T))}}{\|g^\dagger \|_{L^2(\partial\Omega\times(0,T))}},\\
     &e_{u_0} =\frac{\|u_0^\dagger-u_0\|_{L^2(\Omega)}}{\|u_0^\dagger\|_{L^2(\Omega)}},  &&e_{u}=\frac{\|u^\dagger-u\|_{L^2( \Omega )}}{\|u^\dagger \|_{L^2( \Omega)}},  &&&e_{\mathcal{C}}=\frac{\|\mathcal{C}^\dagger-\mathcal{C}\|_{L^2( \Omega\times(0,T))}}{\|\mathcal{C}^\dagger \|_{L^2( \Omega\times(0,T))}},
\end{align*}
where we use $\dagger$ to denote the ground truth.

\begin{figure}
     \centering	
     \begin{tabular}{cccc}
		\includegraphics[width=0.20\textwidth]{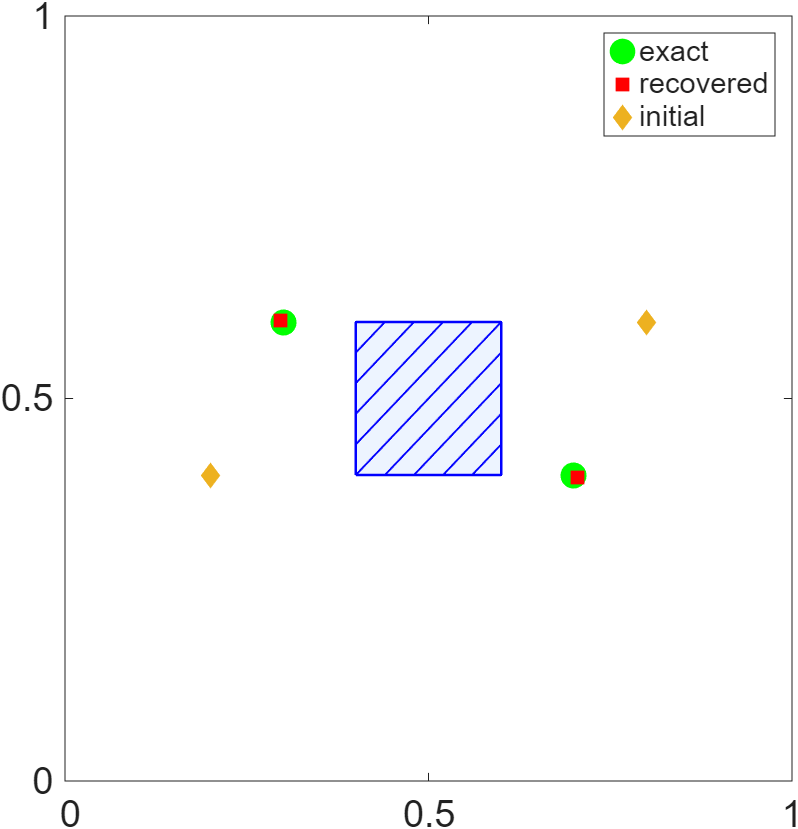}&
		\includegraphics[width=0.24\textwidth]{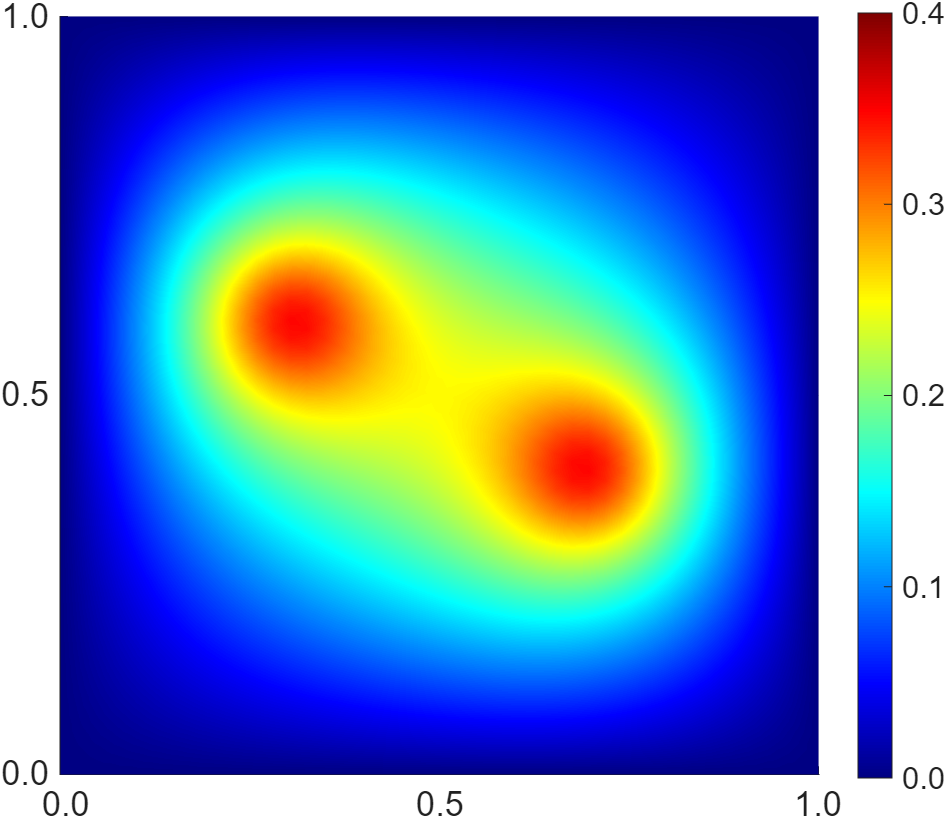}  &
		\includegraphics[width=0.24\textwidth]{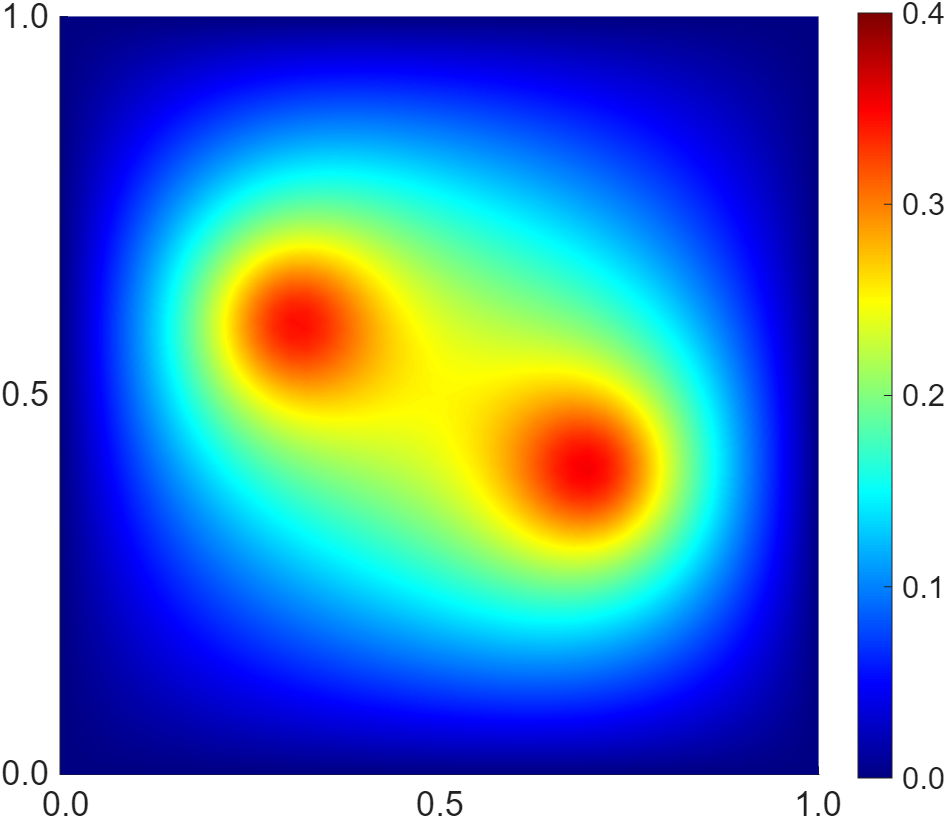} &
        \includegraphics[width=0.25\textwidth]{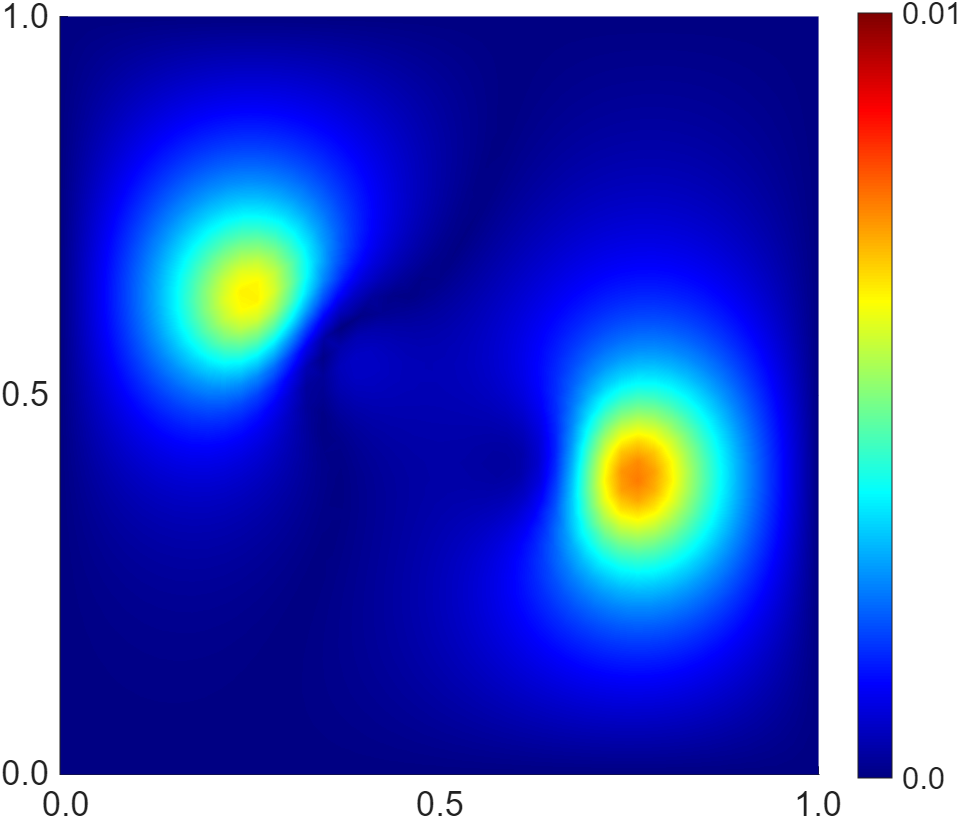}\\
        \includegraphics[width=0.20\textwidth]{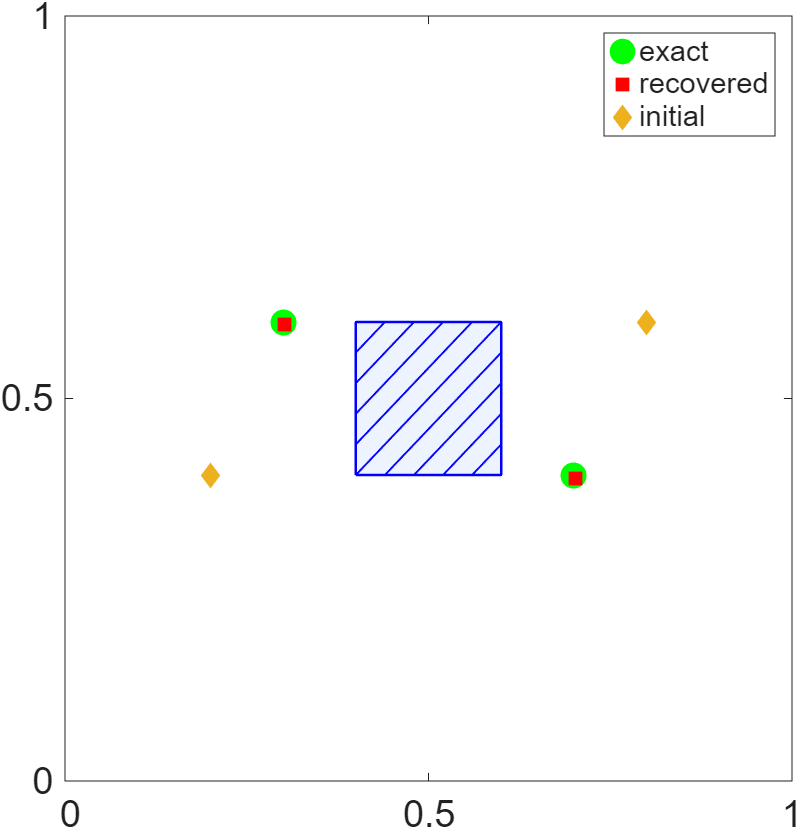}&
		\includegraphics[width=0.24\textwidth]{pics/example_adv_ue_B0.1.png}  &
		\includegraphics[width=0.24\textwidth]{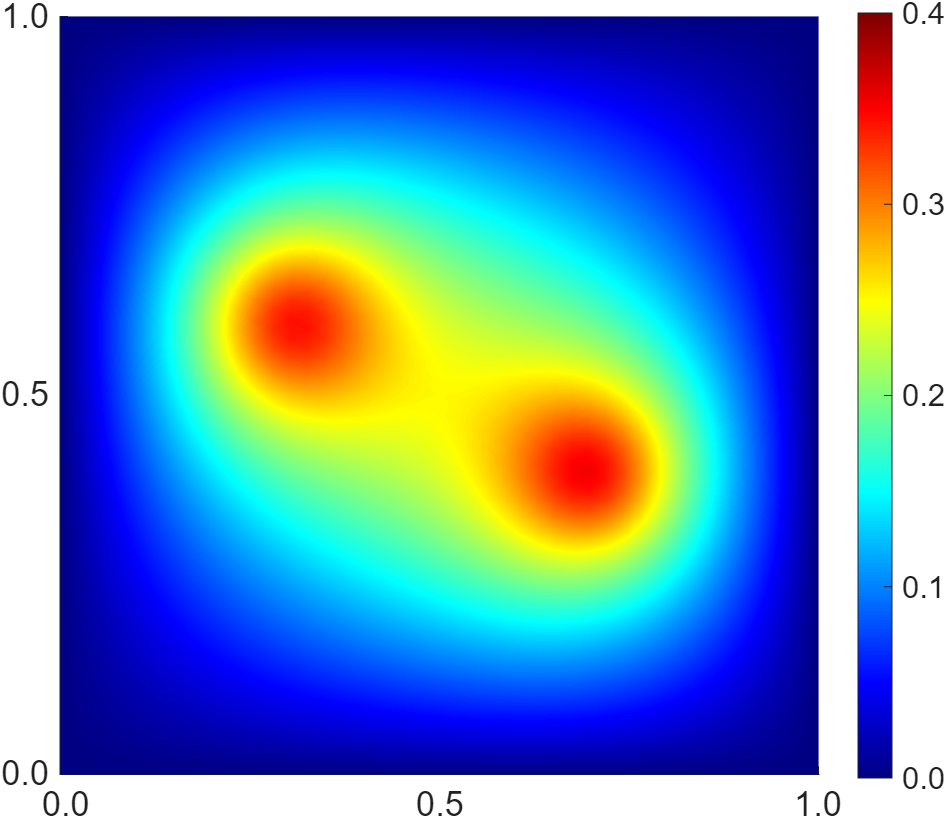} &
        \includegraphics[width=0.25\textwidth]{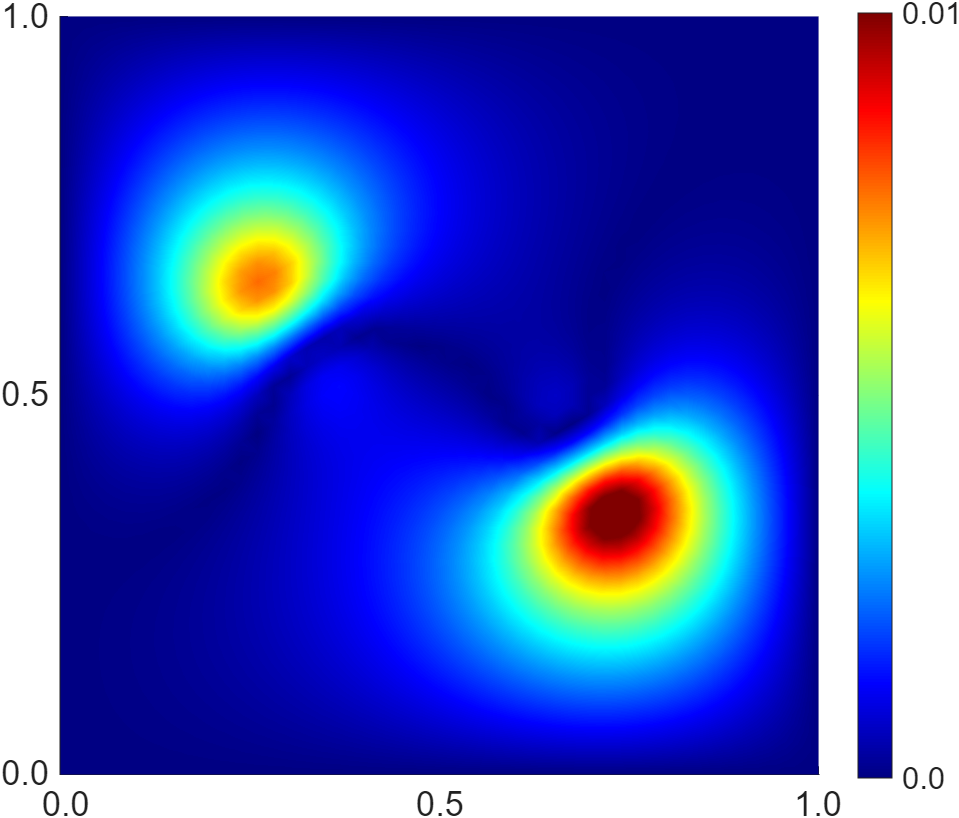}\\
        (a) source locations & (b) true $u^\dagger$ & (c) reconstructed $u$ & (d)  error $|u^\dagger-u|$
	\end{tabular}  
    \caption{The reconstructions  for Example~\ref{ex:2d_elliptic} at different noise levels. First row: $1\%$ noise. Second row: $5\%$ noise.}	\label{Fig:2d_elliptic}
 \end{figure}

To improve the efficiency of the numerical reconstruction, we assume the boundary data $g$ and the initial condition $u_0$ is known in Example~\ref{ex:2d_elliptic}--\ref{ex:3d_elliptic}. 
We start with a simple example in the steady-state case \eqref{equabv}.
\begin{example}\label{ex:2d_elliptic}
Let $d(x_1,x_2)=1+x_1 (1-x_1)\sin(\pi x_2) $, $g(x_1,x_2)=0$, $V=(V_1,V_2)^T$ with
\begin{itemize}
    \item[(i)] $V_1(x_1,x_2)=0.4x_2(1-x_2)$ and  $V_2(x_1,x_2)=0$,
    \item[(ii)]$V_1(x_1,x_2)=40x_2(1-x_2)$ and  $V_2(x_1,x_2)=0$.
\end{itemize}
We aim to reconstruct the unknown pollution source $\lambda_1\delta_{S_1}+\lambda_2\delta_{S_2}$ and the corresponding pollutant concentration field $u$. The ground truths are $\lambda_1=\lambda_2=1$, $S_1=(0.3,0.6)$ and $S_2=(0.7,0.4)$. The measurement domain is $\omega=(0.4,0.6)^2$.
\end{example}
In this example, we parameterize the coordinates of the two point sources as $\bm{\mu}$ and set $\bm{\lambda}=(\lambda_1,\lambda_2)$. The penalty parameter is chosen as $\lambda_{\bm{\mu}}^k\equiv \text{1e-3}$. We first consider case (i), in which the velocity field $V$ is relatively small compared with the diffusion coefficient $b$, leading to a diffusion-dominated regime. The reconstruction results are presented in Figure~\ref{Fig:2d_elliptic} and Table~\ref{tab:error}. In Figure~\ref{Fig:2d_elliptic}(a), we show the reconstructed pollution source. The reconstruction is accurate and robust with respect to noise, even when the initial guess is far from the ground truth. In Figure~\ref{Fig:2d_elliptic}(b)--(d), we show the reconstructed pollutant concentration field. The algorithm recovers the field over the whole domain from measurements in the small domain $\omega$, demonstrating the effectiveness of the proposed method. 

\begin{figure}
     \centering	
     \begin{tabular}{cccc}
		\includegraphics[width=0.20\textwidth]{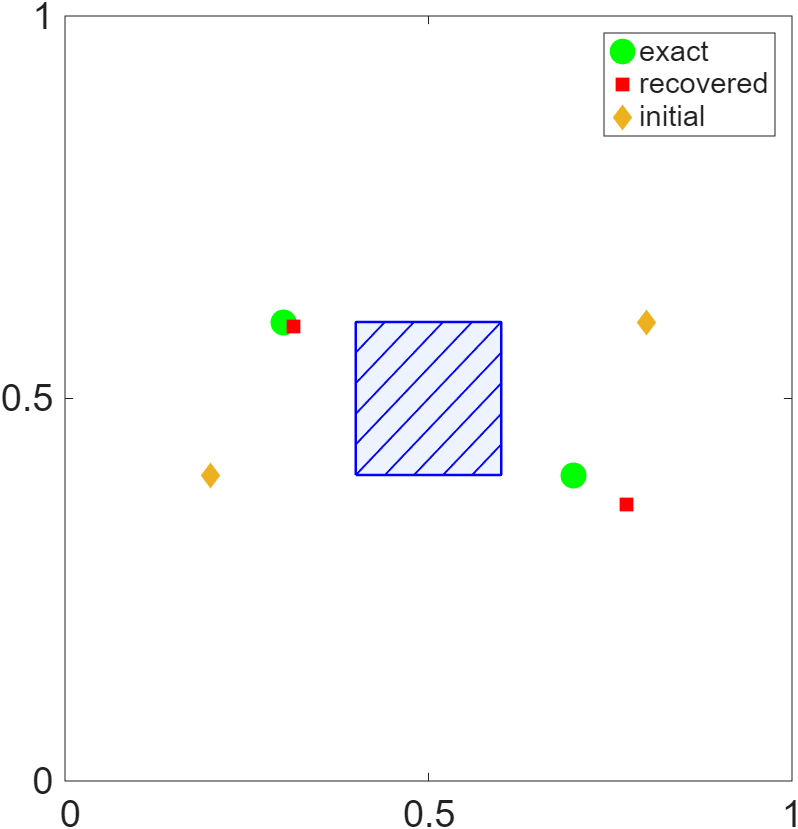}&
		\includegraphics[width=0.24\textwidth]{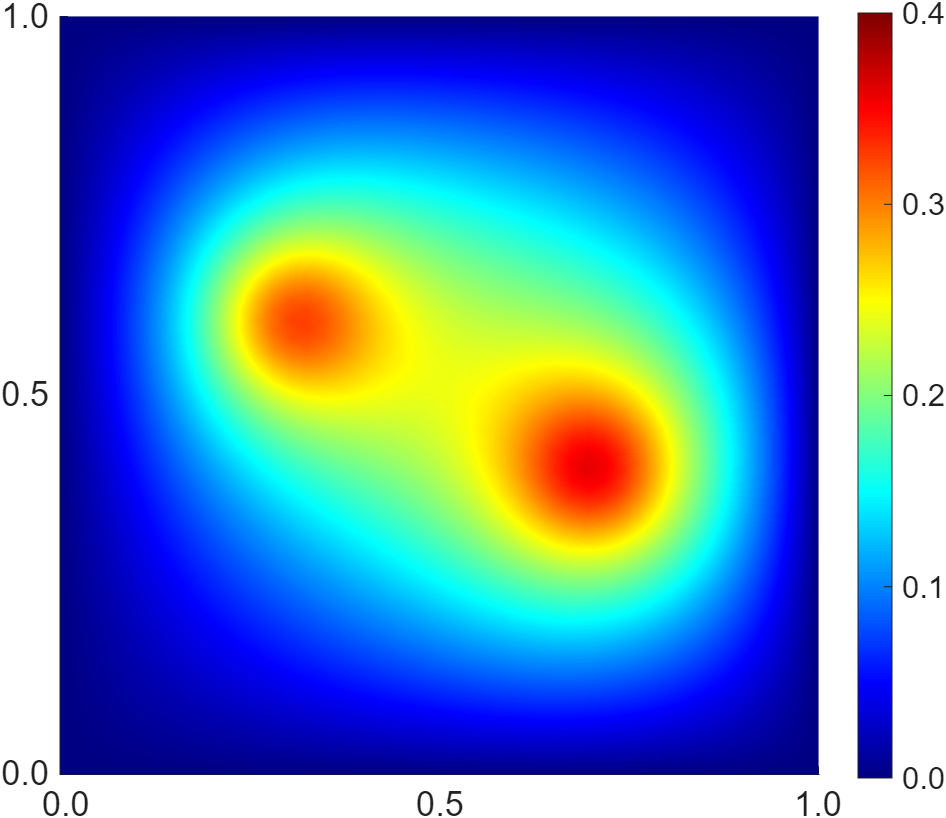}  &
		\includegraphics[width=0.24\textwidth]{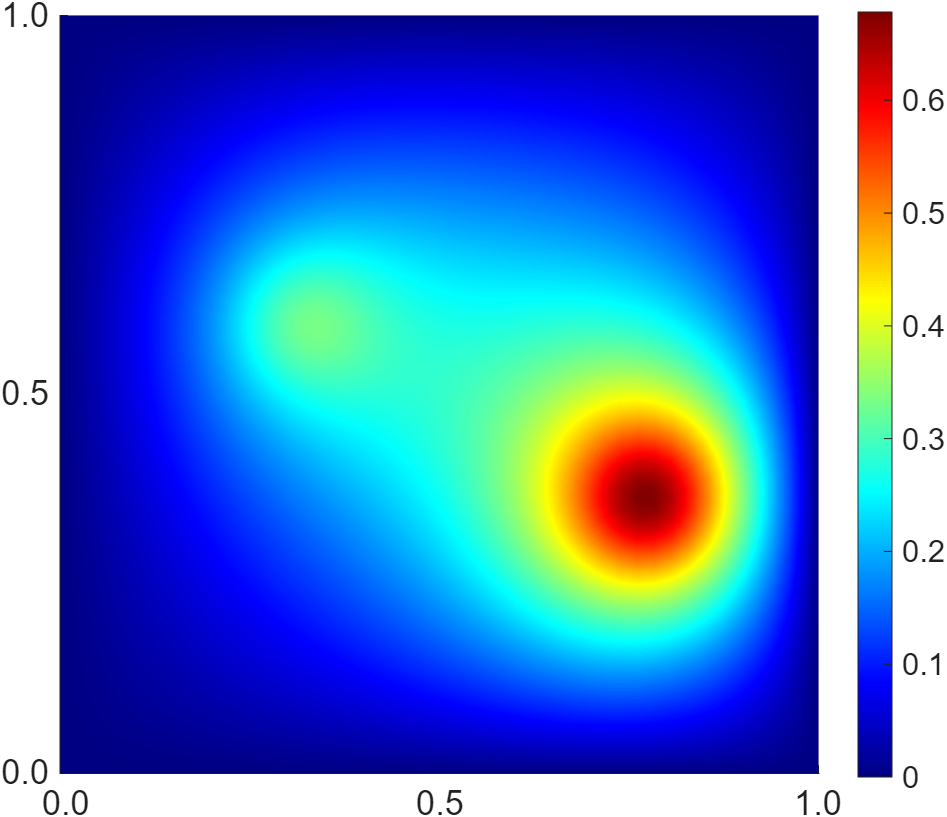} &
        \includegraphics[width=0.24\textwidth]{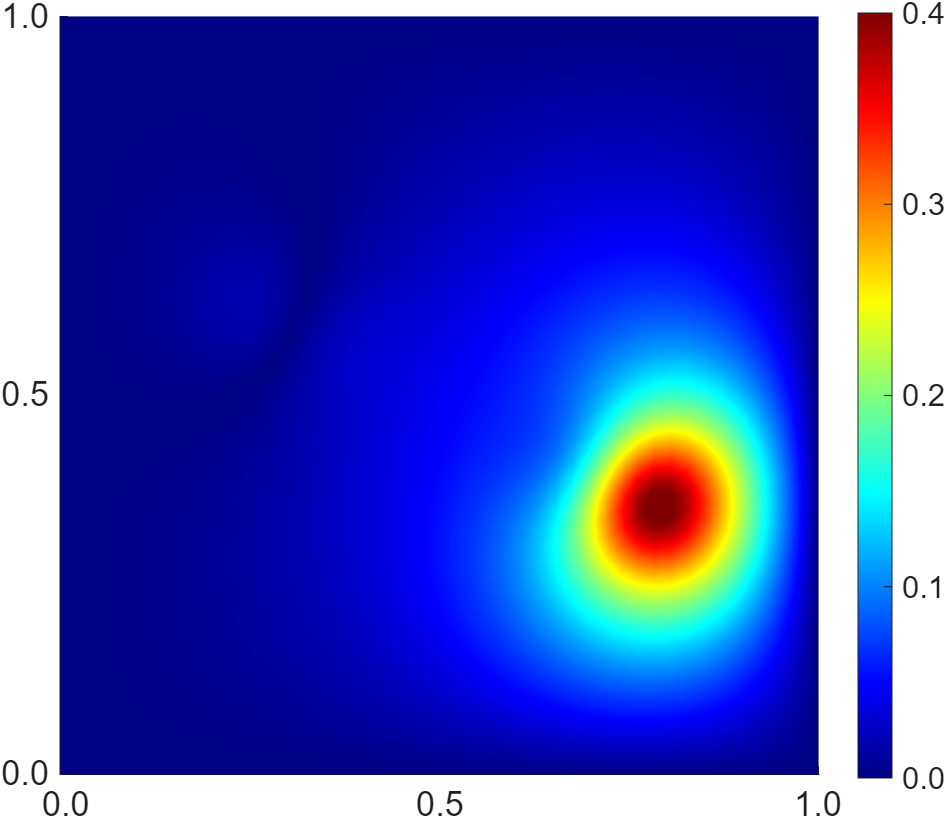} \\
        (a) source locations & (b) true $u^\dagger$ & (c) reconstructed $u$ & (d)  error $|u^\dagger-u|$
	\end{tabular}  
    \caption{The reconstructions  for Example~\ref{ex:2d_elliptic} in $1\%$ noise at advection dominated case. }	\label{Fig:2d_elliptic_adv}
 \end{figure}

We then consider case (ii), an advection-dominated regime with a relatively large velocity field. Here, we employ the streamline upwind/Petrov–Galerkin (SUPG) formulation \cite{Hughes:1982} to stabilize the advection term in the proposed discretization. The reconstruction results are presented in Figure~\ref{Fig:2d_elliptic_adv} and Table~\ref{tab:error}. In this high-advection setting, the reconstruction exhibits a bias along the velocity field: the source $S_1$ is recovered relatively well, whereas the reconstruction of $S_2$ is severely deteriorated. This asymmetry is directional: the observation window $\omega$ lies downstream of $S_1$, so the flow carries its signature into the data, but upstream of $S_2$, whose signature can reach the data only against the flow.
This numerical instability is consistent with the theoretical analysis: the unique continuation from $\omega$ becomes unstable in the advection-dominated regime, the constants in the stability estimates growing exponentially with $|V|$.

 \begin{table}[hbt!]
  \centering
  \begin{threeparttable}
    \caption{The relative errors for Example~\ref{ex:2d_elliptic}-\ref{ex:bdry} at different noise levels.\label{tab:error}}
    \begin{tabular}{l|cccc|cccc}
        \toprule
        & \multicolumn{4}{c}{$\delta=1\%$} & \multicolumn{4}{c}{$\delta=5\%$}\\
\cmidrule{2-5} \cmidrule{6-9}
           &   $e_{\bm{\mu}}$ & $e_{\bm{\lambda}}$ & $e_{g}$ & $e_{u}(e_{\mathcal{C}})$ & $e_{\bm{\mu}}$ & $e_{\bm{\lambda}}$ & $e_g$ & $e_{u}(e_{\mathcal{C}})$\\
        \midrule
           Example~\ref{ex:2d_elliptic}(i)  & 3.95e-3 & 1.70e-2 & / & 1.75e-3  & 5.67e-3 & 2.06e-2 & / & 2.98e-3 \\
           Example~\ref{ex:2d_elliptic}(ii) & 7.89e-2 & 9.31e-1 & / & 3.82e-1 & 7.91e-2 & 9.21e-1 & / & 3.78e-1\\
           Example~\ref{ex:2d_para}   &  3.18e-3 & 7.72e-3 & / & 4.93e-3& 5.40e-3 & 1.09e-2 & / & 5.96e-3 \\
           Example~\ref{ex:3d_elliptic}  & 8.25e-3 & 3.20e-2 & / & 1.80e-2& 3.56e-2 & 1.44e-1 & / & 1.07e-1  \\
           Example~\ref{ex:bdry}  &   3.99e-2 & 2.83e-1 & 3.86e-1& 7.10e-2 & 4.11e-2 & 2.93e-1 & 6.32e-1& 7.99e-2  \\ 
        \bottomrule
    \end{tabular}
      \end{threeparttable}
\end{table}

\begin{figure}
     \centering	
     \begin{tabular}{cccc}
		\includegraphics[width=0.21\textwidth]{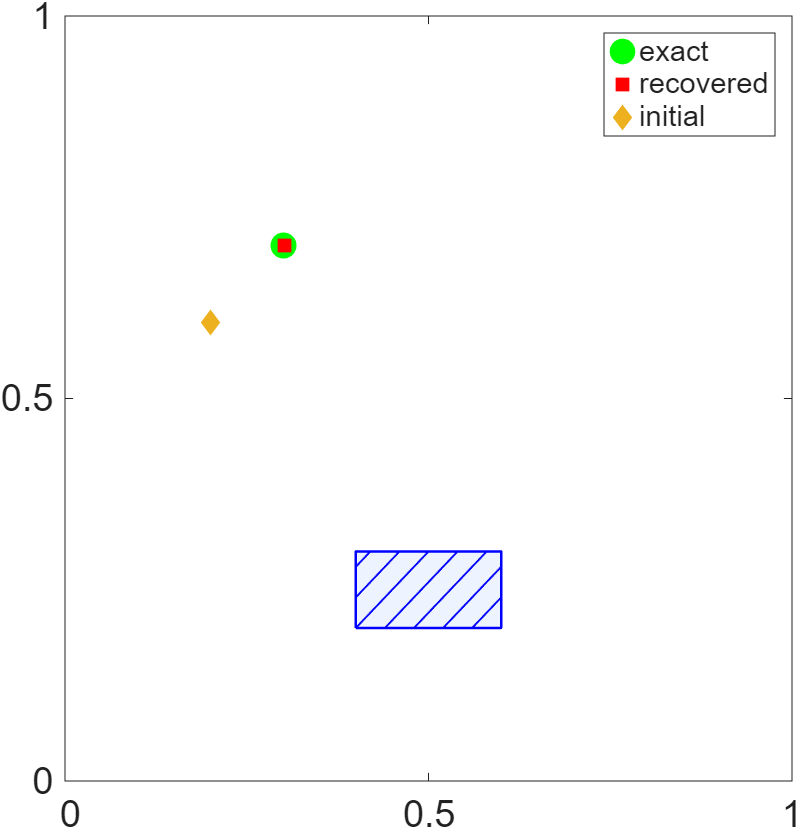}&
		\includegraphics[width=0.21\textwidth]{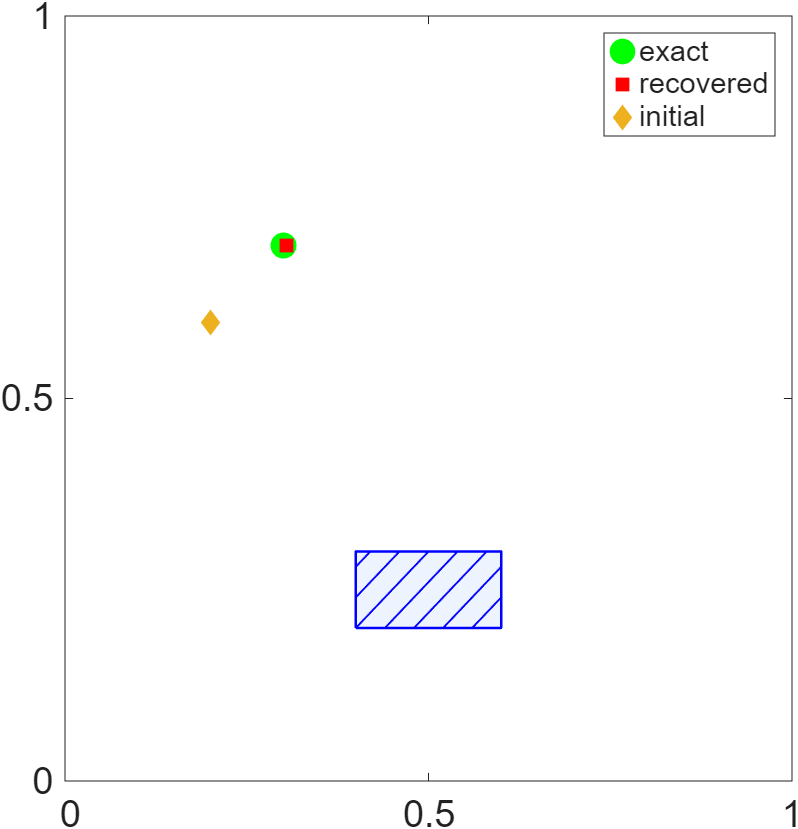}  &
		\includegraphics[width=0.24\textwidth]{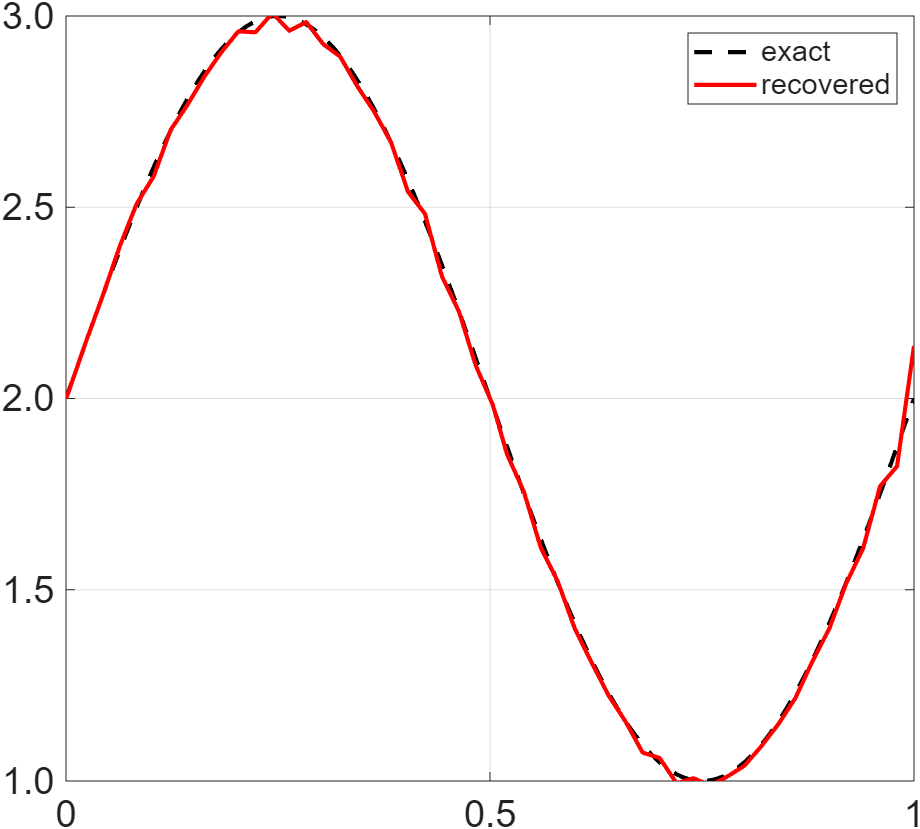} &
        \includegraphics[width=0.24\textwidth]{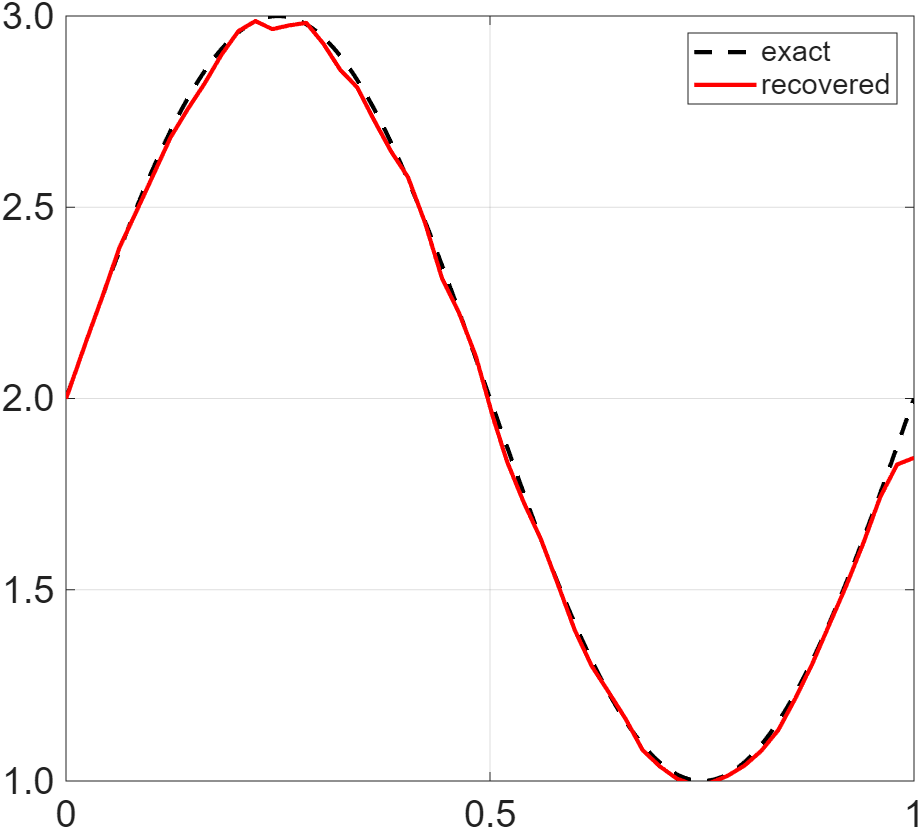}\\ 
        (a) $S_1$ with $\delta=1\%$  & (b) $S_1$ with $\delta=5\%$  & (c) $\lambda_1(t)$ with $\delta=1\%$  & (d)   $\lambda_1(t)$ with $\delta=5\%$
	\end{tabular}  
    \caption{The reconstructions of pollution sources  for Example~\ref{ex:2d_para} at different noise levels.  }	\label{Fig:2d_para_source}
 \end{figure}

\begin{figure}
     \centering	
     \begin{tabular}{ccc} 
		\includegraphics[width=0.24\textwidth]{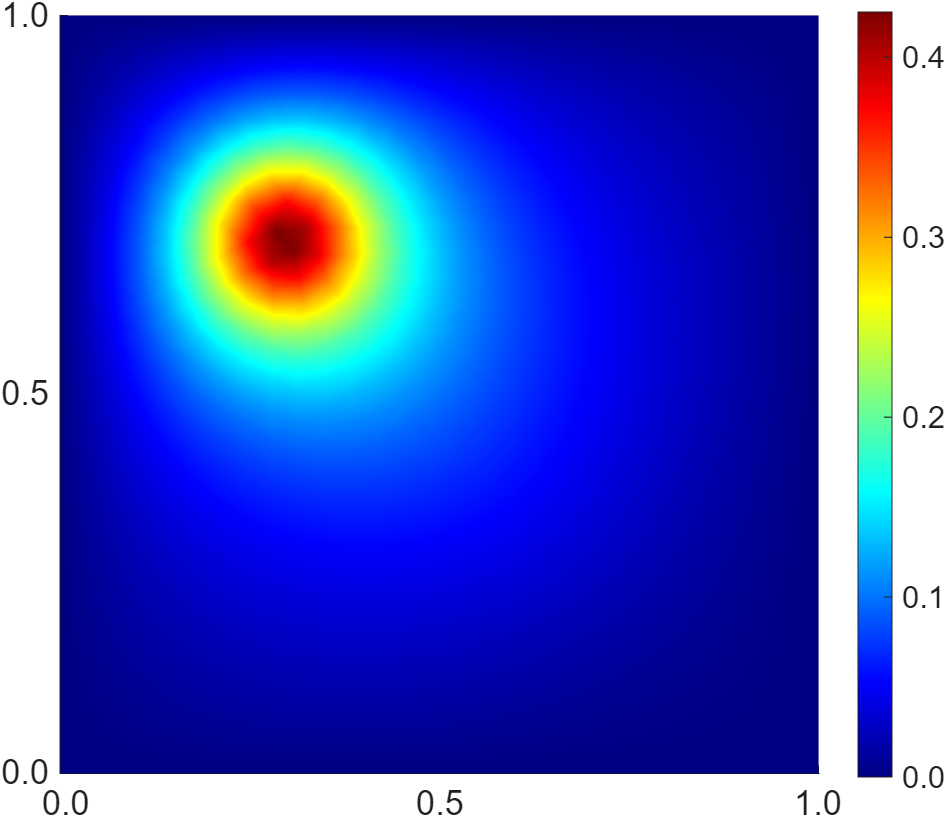}  &
		\includegraphics[width=0.24\textwidth]{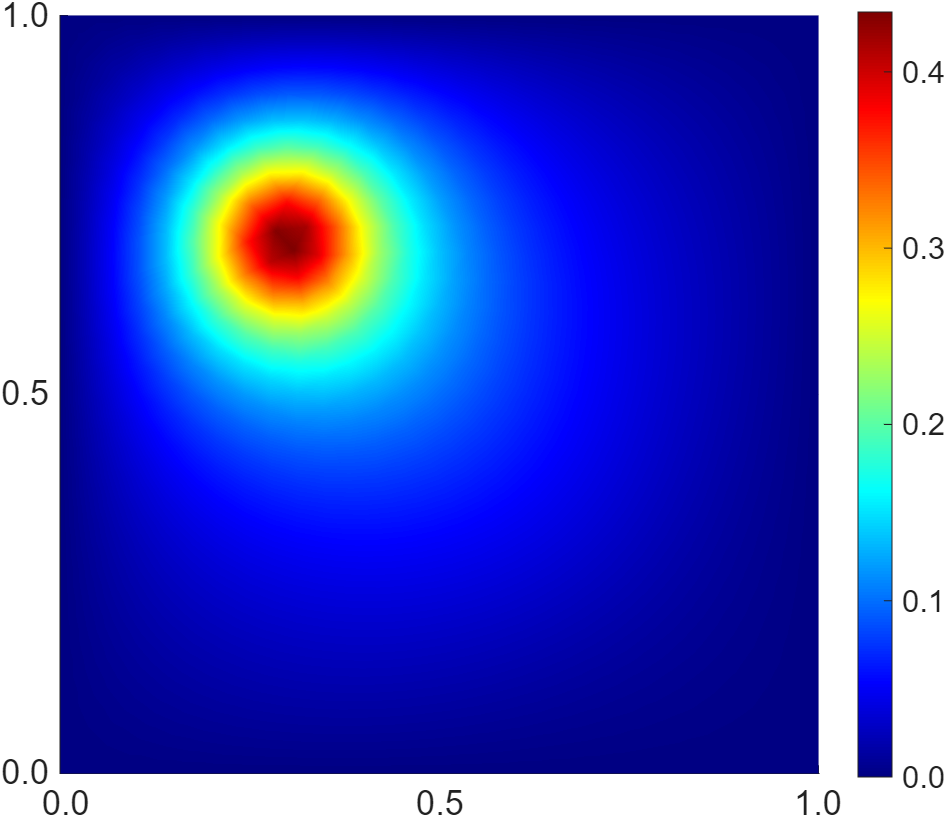} &
        \includegraphics[width=0.24\textwidth]{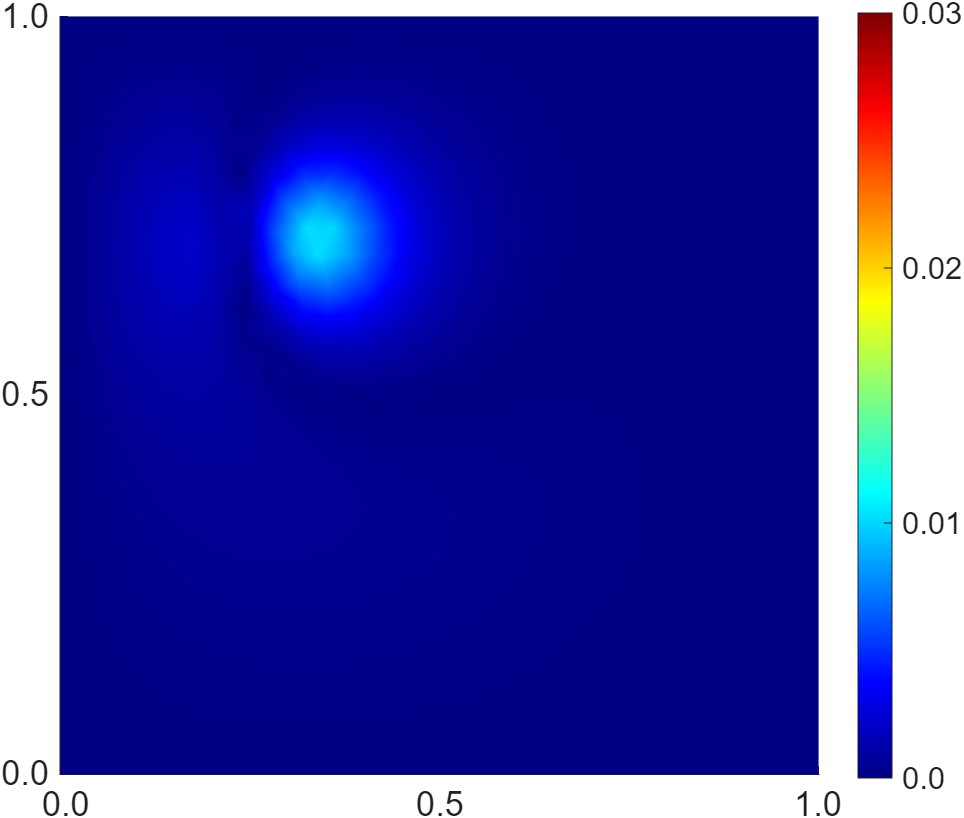} \\
        \includegraphics[width=0.24\textwidth]{pics/example_para_ue.png}  &
		\includegraphics[width=0.24\textwidth]{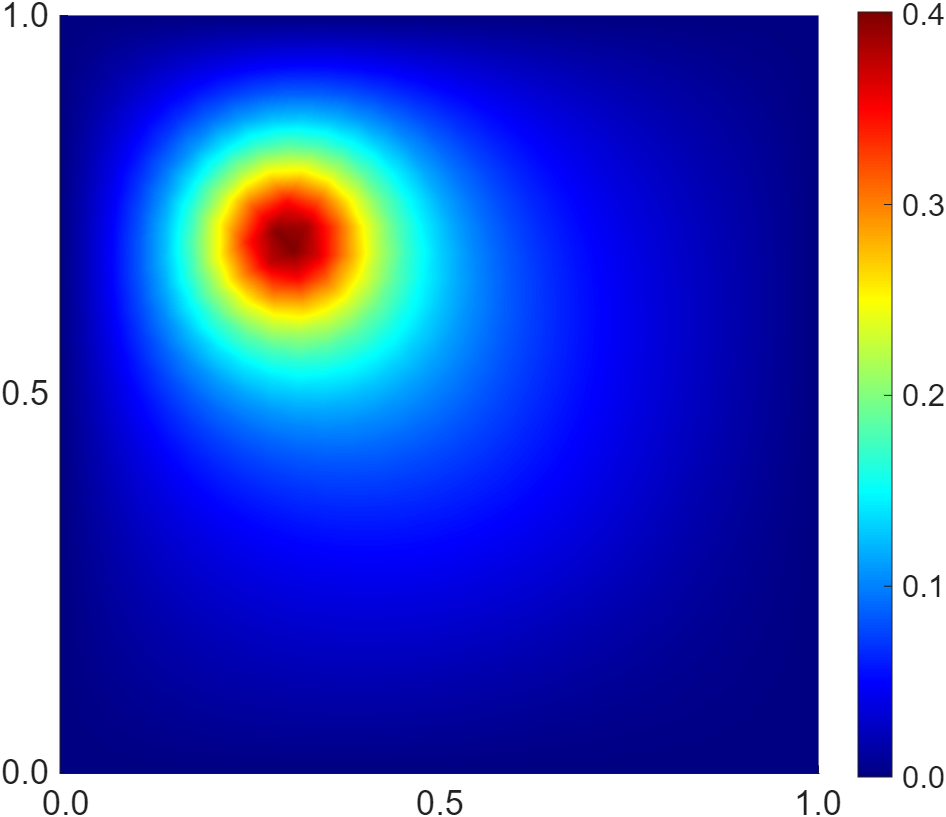} &
        \includegraphics[width=0.24\textwidth]{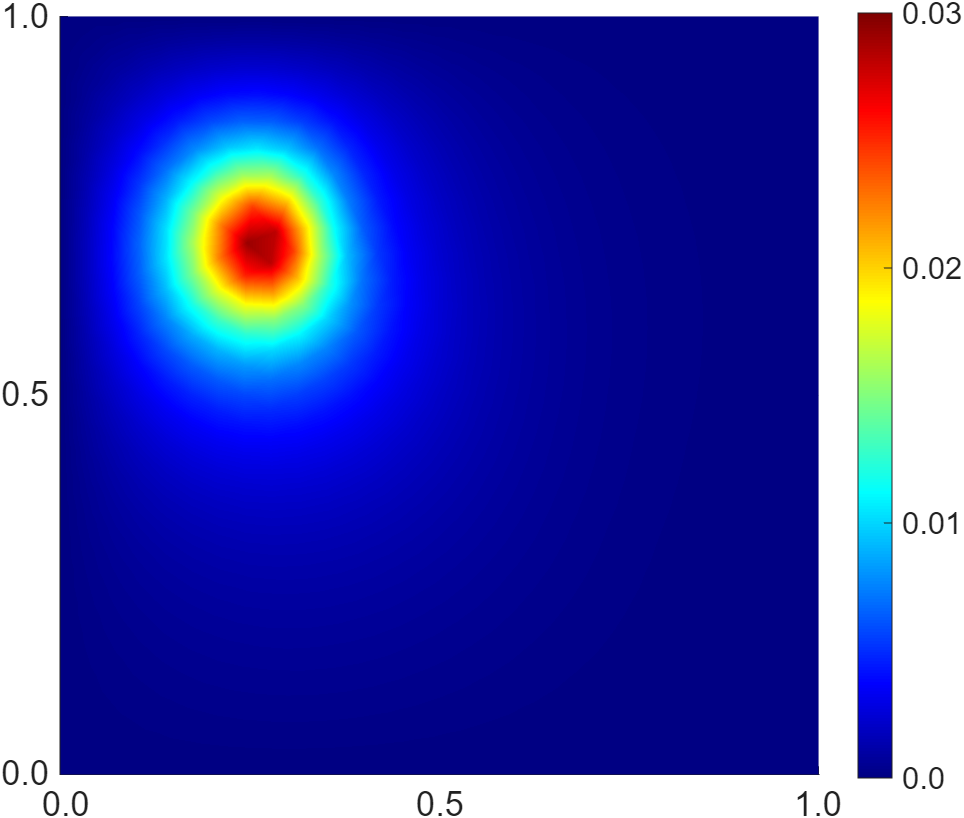} \\
         (a) true $\mathcal{C}^\dagger$ & (b) reconstructed $\mathcal{C}$ & (c)  error $|\mathcal{C}^\dagger-\mathcal{C}|$
	\end{tabular}  
    \caption{The reconstructions of pollutant concentration fields $\mathcal{C}$ at $T=1$ for Example~\ref{ex:2d_para} at different noise levels. First row: $1\%$ noise. Second row: $5\%$ noise.}	\label{Fig:2d_para_sol}
 \end{figure}

Next we consider a time-dependent problem \eqref{equa}.
\begin{example}\label{ex:2d_para}
Let $d(x_1,x_2)=2+ \sin(2\pi x_1)\sin(2\pi x_2) $, $V=(V_1,V_2)^T$ with $V_1(x_1,x_2)=0.4x_2(1-x_2)$ and  $V_2(x_1,x_2)=0$, $g(x_1,x_2,t)=0$ and $u_0(x_1,x_2)=\frac{1}{2}\sin(2\pi x_1)\sin(2\pi x_2)$. We aim to reconstruct the unknown pollution source $\lambda_1(t)\delta_{S_1} $ and the corresponding pollutant concentration field $\mathcal{C}$. The ground truths are $\lambda_1(t)=2+\sin(2\pi t)$ and $S_1=(0.3,0.7)$. The measurement domain is $\omega=(0.4,0.6)\times(0.2,0.3)$.
\end{example}
Here the location of the point source is parameterized by $\bm{\mu}$. The time-dependent intensity $\lambda_1(t)$ is approximated by a continuous piecewise linear function on $[0,T]$, $\lambda_1(t)\approx\sum_{i=1}^{N_T} \lambda_1^i\,\phi_i(t)$, with $N_T=50$ in the experiments. The penalty parameter is set to $\lambda_{\bm{\mu}}^k\equiv \text{1e-6}$. 
The reconstruction results are presented in Figures~\ref{Fig:2d_para_source}, \ref{Fig:2d_para_sol} and Table~\ref{tab:error}. Both the time-dependent pollution source and the pollutant concentration field are recovered accurately even at high noise levels, which illustrates the robustness of the proposed method.

\begin{figure}
     \centering	
     \begin{tabular}{cccc}
		\includegraphics[width=0.21\textwidth]{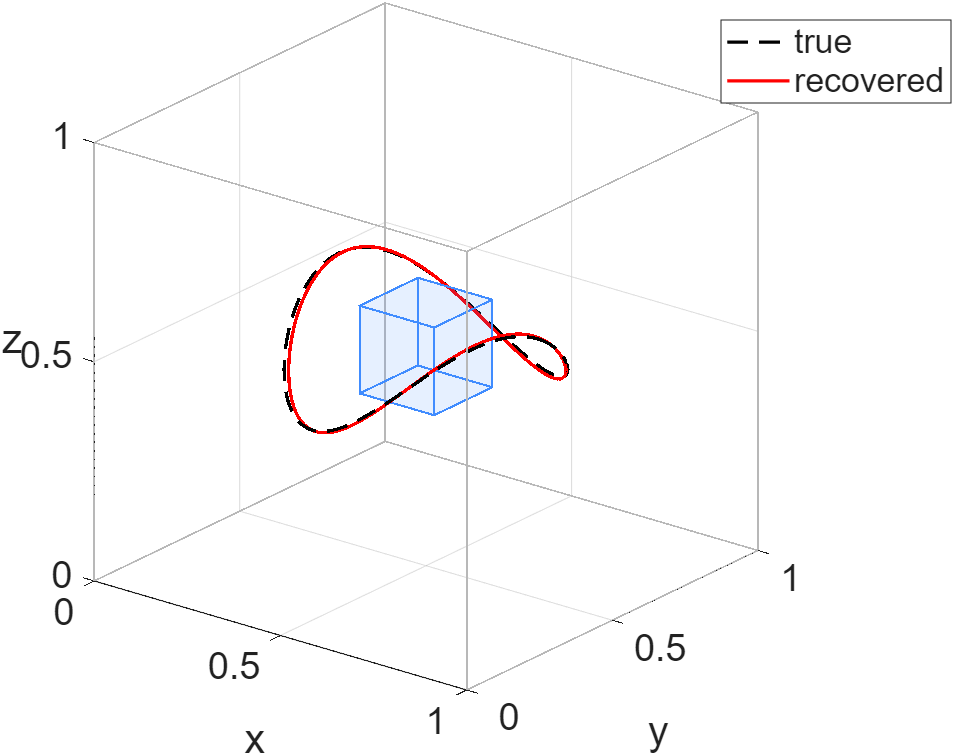}&
		\includegraphics[width=0.21\textwidth]{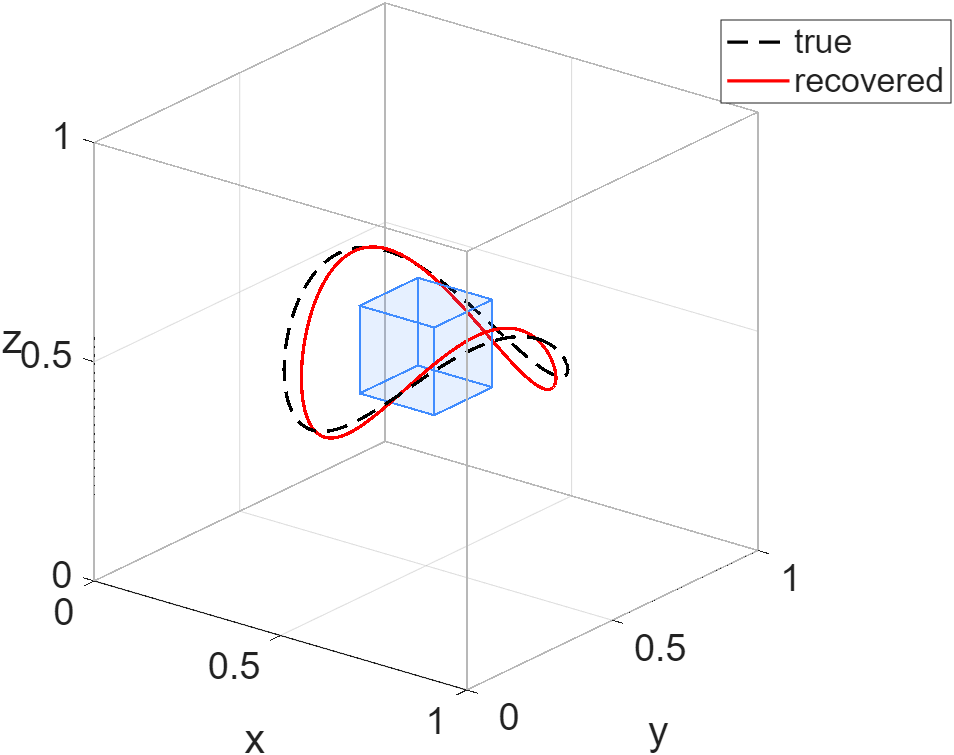}  &
		\includegraphics[width=0.21\textwidth]{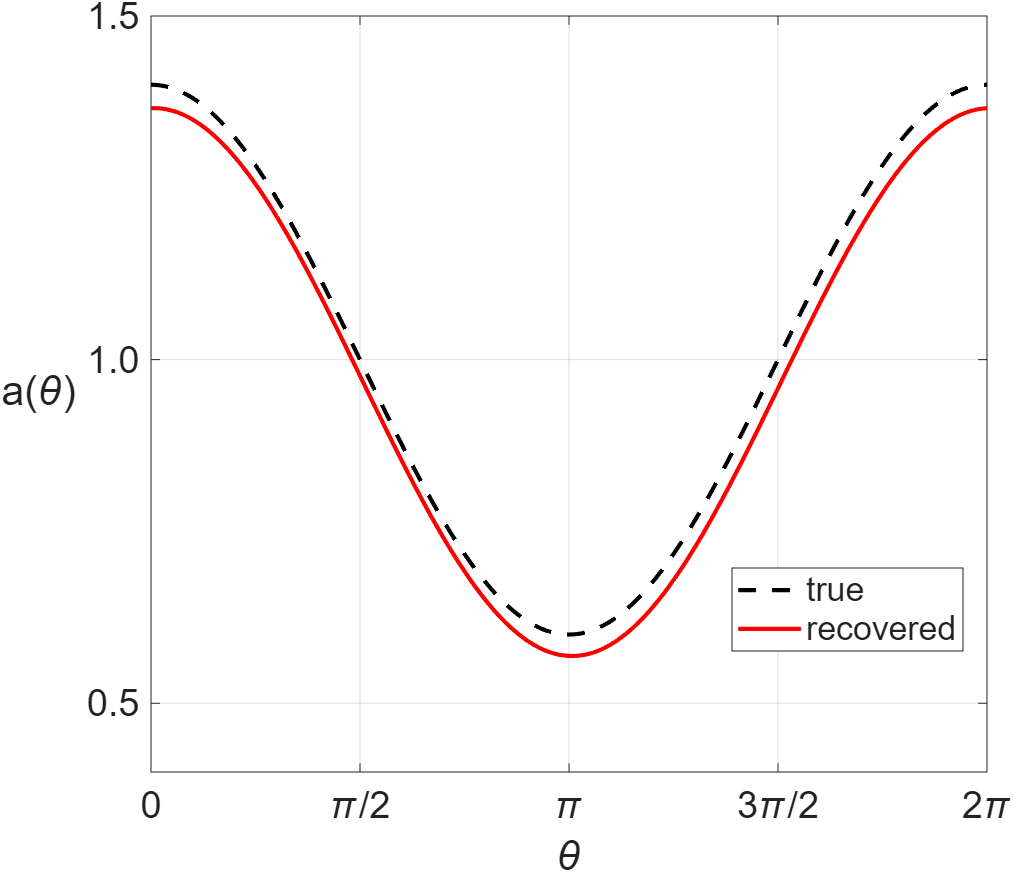} &
        \includegraphics[width=0.21\textwidth]{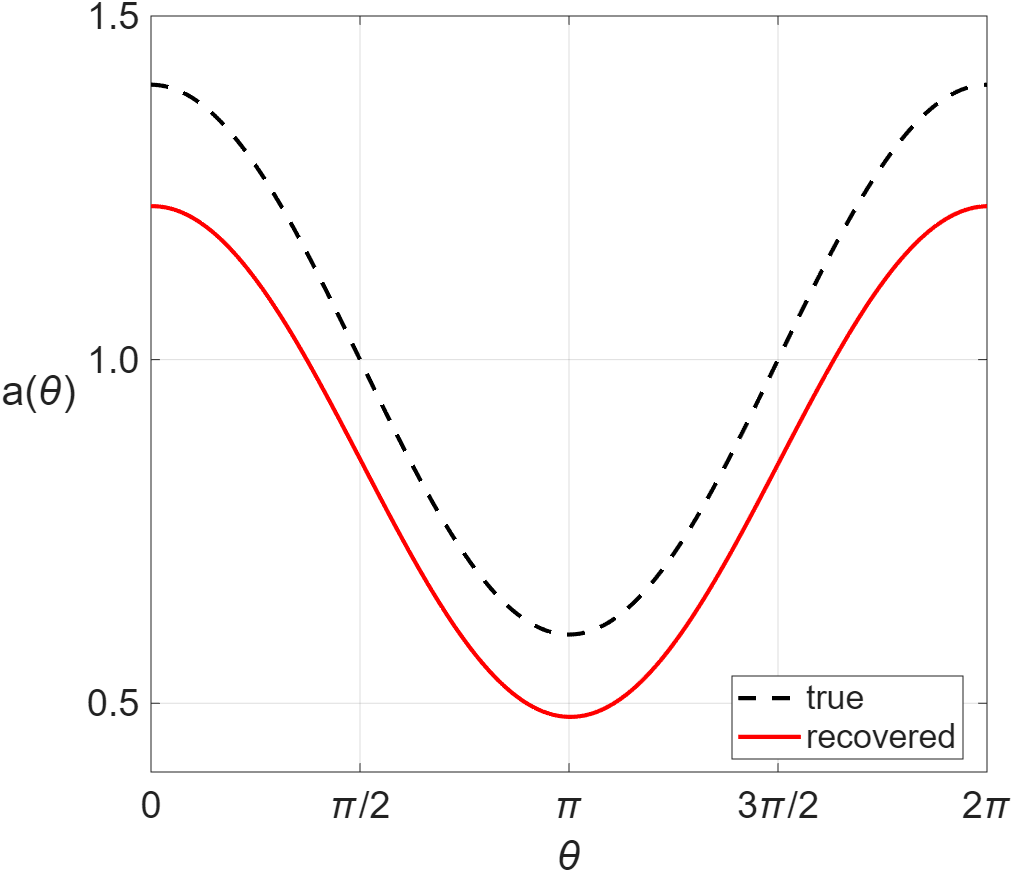}\\ 
        (a) $\gamma$ with $\delta=1\%$  & (b) $\gamma$ with $\delta=5\%$  & (c) $\mu$ with $\delta=1\%$  & (d)   $\mu$ with $\delta=5\%$
	\end{tabular}  
    \caption{The reconstructions of pollution sources  for Example~\ref{ex:3d_elliptic} at different noise levels.  }	\label{Fig:3d_elliptic_source}
 \end{figure}

\begin{figure}
     \centering	
     \begin{tabular}{ccc} 
		\includegraphics[width=0.24\textwidth]{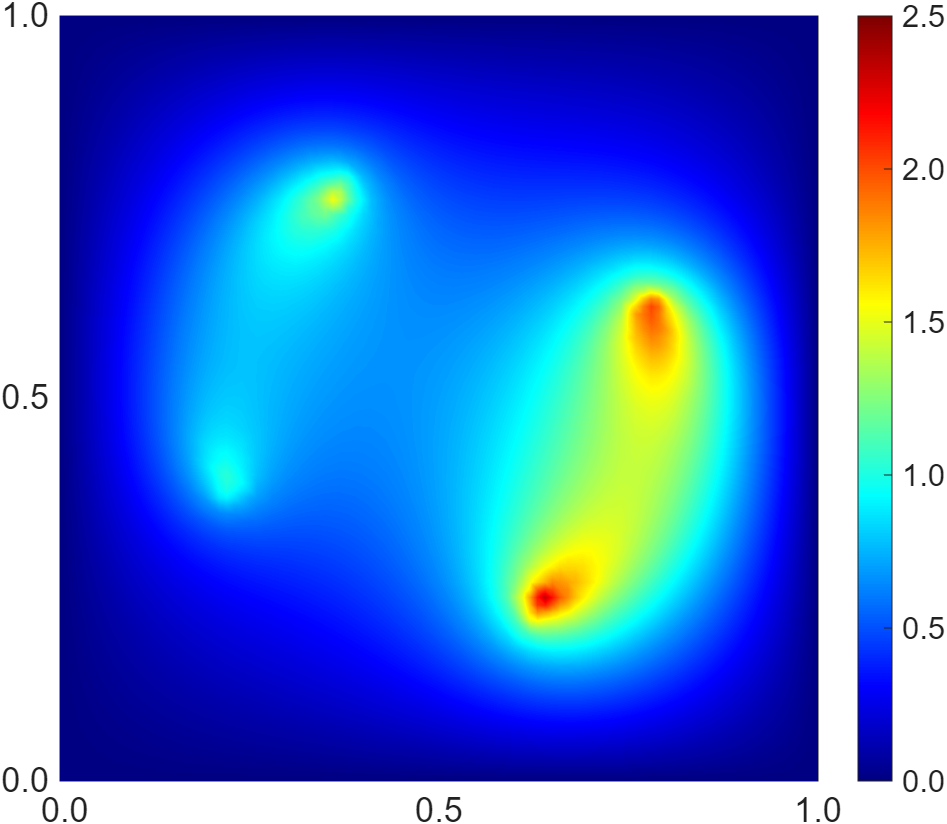}  &
		\includegraphics[width=0.24\textwidth]{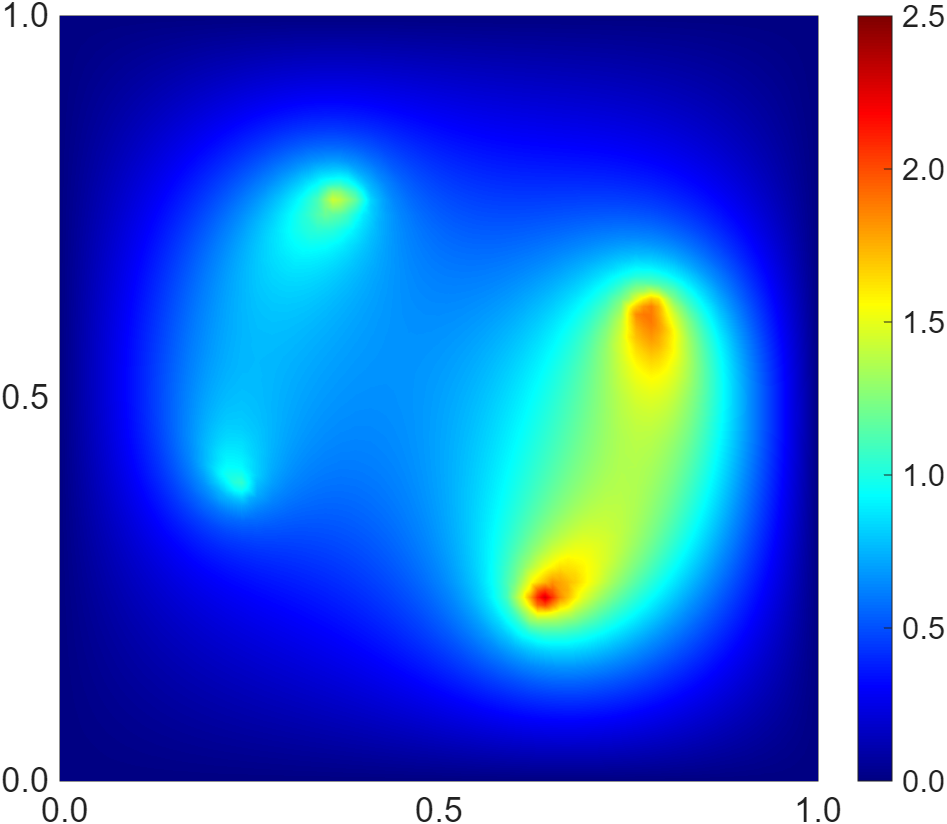} &
        \includegraphics[width=0.24\textwidth]{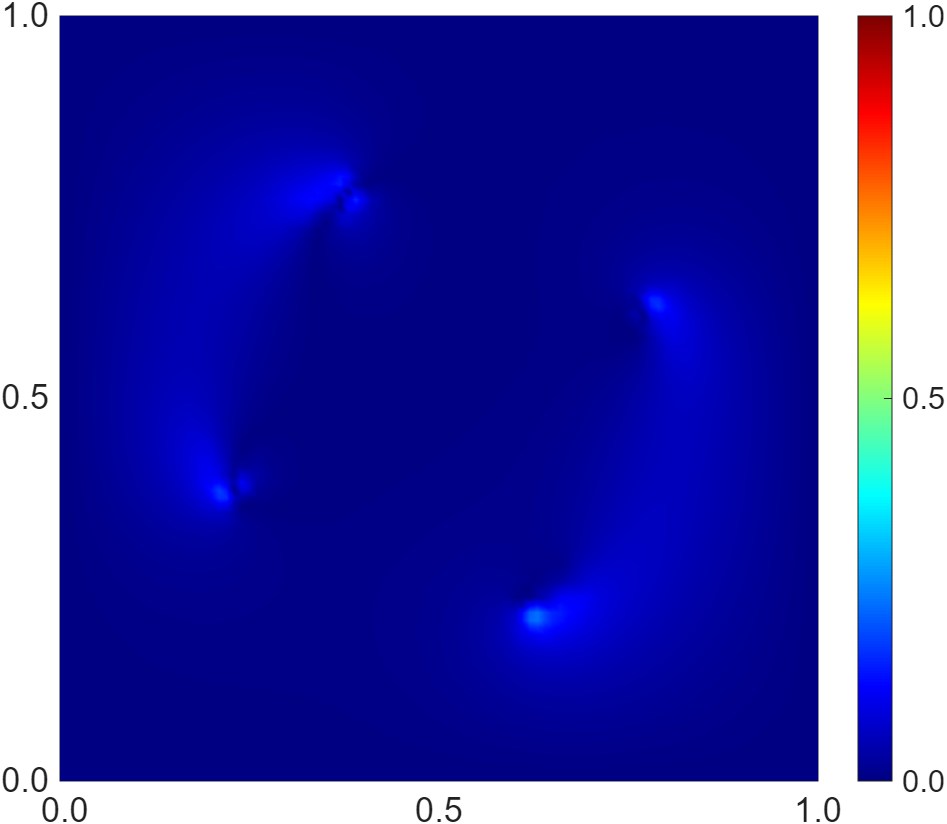} \\
        \includegraphics[width=0.24\textwidth]{pics/example_3d_ue.png}  &
		\includegraphics[width=0.24\textwidth]{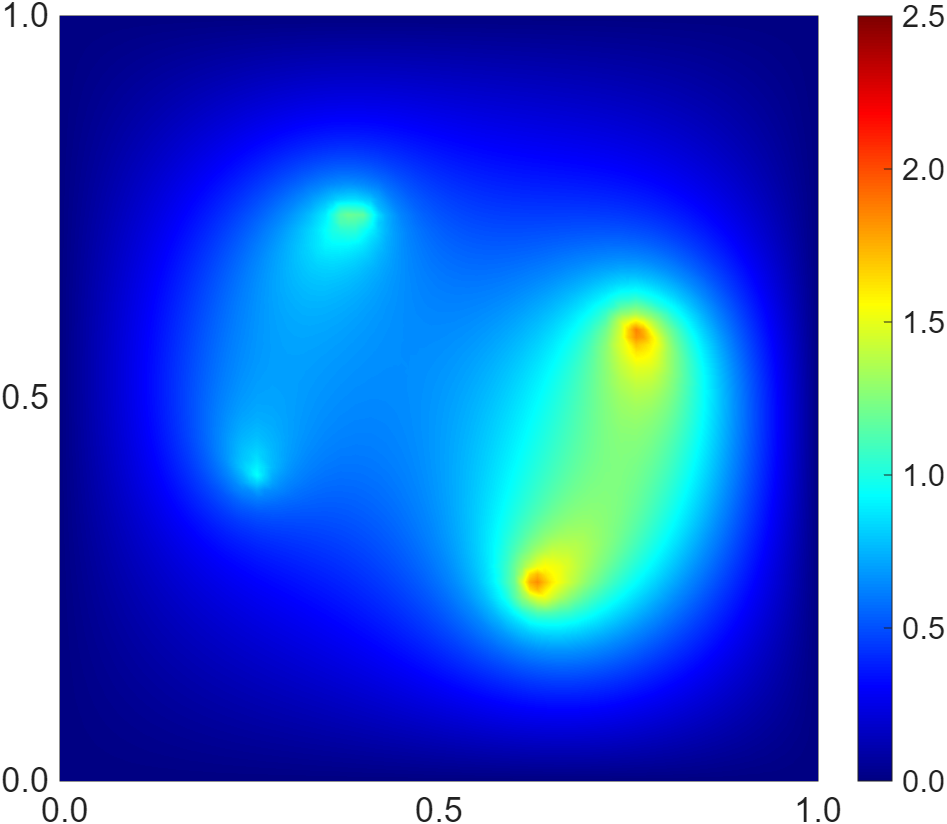} &
        \includegraphics[width=0.24\textwidth]{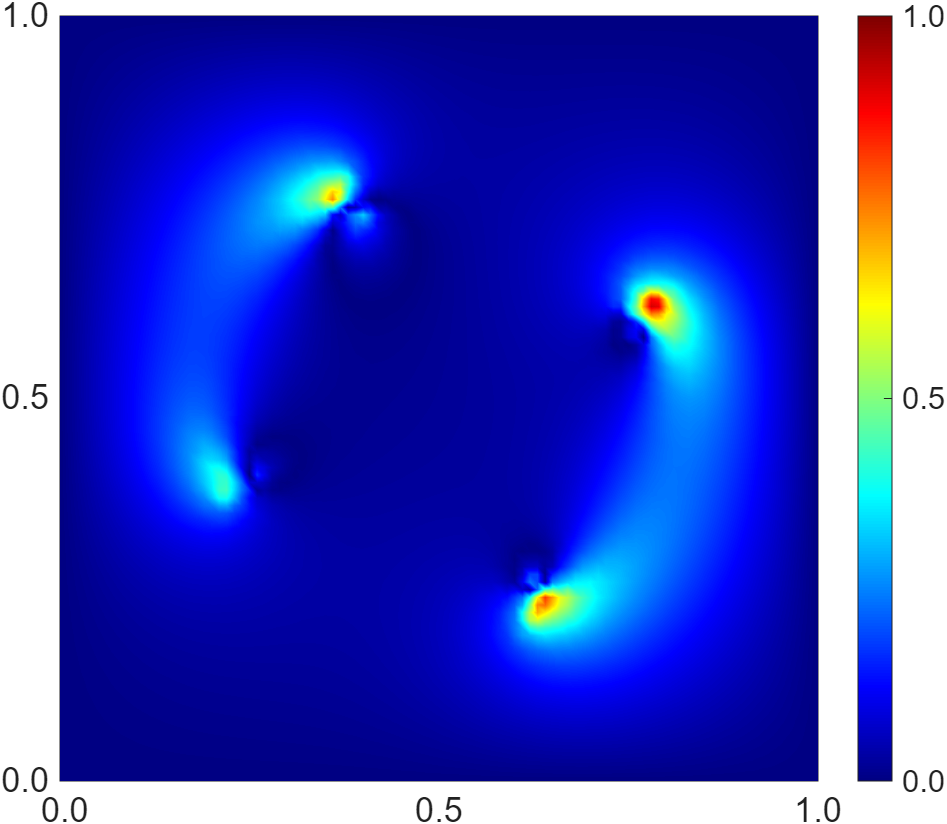} \\
         (a) true $u^\dagger$ & (b) reconstructed $u$ & (c)  error $|u^\dagger-u|$
	\end{tabular}  
    \caption{The reconstructions of pollutant concentration fields $u$ at the slice $x_3=0.5$ for Example~\ref{ex:3d_elliptic} at different noise levels. First row: $1\%$ noise. Second row: $5\%$ noise.}	\label{Fig:3d_elliptic_sol}
 \end{figure}

We now turn to a three-dimensional steady-state problem \eqref{equabv}, in which the pollution source is supported on a curve. For this example, the data are generated on a fine grid with mesh size $h=\frac{1}{40}$, and the inversion is solved on a coarse grid with $h=\frac{1}{20}$.
\begin{example}\label{ex:3d_elliptic}
Let $d(x_1,x_2,x_3)=1+ 0.3\sin(\pi x_1)\sin(\pi x_2)\sin(\pi x_3) $, $V=(V_1,V_2,V_3)^T$ with $V_1(x_1,x_2,x_3)=0.4x_2(1-x_2)$ and  $V_2(x_1,x_2,x_3)=V_3(x_1,x_2,x_3)=0$, $g(x_1,x_2,x_3)=0$. We aim to reconstruct the unknown pollution source $ \mu $ and the corresponding pollutant concentration field $u$. The ground truth of $\mu$ is supported on the curve $\gamma(\theta)=  ( 0.5 + 0.3\cos\theta,  0.5  + 0.3 \sin\theta,  0.5 + 0.12\cos(2\theta + 0.7) ) 
$, $\theta\in[0,2\pi]$ with density $\mu(\theta)\equiv 1+0.4\cos(\theta) $.   The measurement domain is $\omega=(0.4,0.6)^3$.
\end{example}
The support of the source is a circle of radius $R=0.3$ centered at $S=(0.5,0.5,0.5)$ in the $(x_1,x_2)$-plane, oscillating in the $x_3$-direction with amplitude $A=0.12$ and phase $P=0.7$. We collect these geometric parameters in $\bm{\mu}=(S,R,A,P)$, which parameterizes the support of the measure. The density along the curve, with respect to $\theta$, is represented by the trigonometric polynomial
$  \mu(\theta)=\lambda_0+\lambda_1\cos\theta+\lambda_2\sin\theta$,
with coefficients $\bm{\lambda}=(\lambda_0,\lambda_1,\lambda_2)$. The penalty parameter are chosen as  $\lambda_{\bm{\mu}}^k\equiv \text{1e-6}$. The reconstruction results are displayed in Figures~\ref{Fig:3d_elliptic_source}, \ref{Fig:3d_elliptic_sol} and Table~\ref{tab:error}. 
Note that in this case the vector $\bm{\mu}$ involves several parameters with  different sensitivities to the data, which poses a numerical challenge. Consequently, the reconstruction is accurate at low noise levels, but deteriorates at high noise levels.

Finally, we consider a more realistic case, in which the pollutant concentration field is reconstructed without any knowledge of the boundary data.
\begin{example}\label{ex:bdry}
Let $d(x_1,x_2)=2+ \sin(2\pi x_1)\sin(2\pi x_2) $, $V=(V_1,V_2)^T$ with $V_1(x_1,x_2)=4x_2(1-x_2)$ and  $V_2(x_1,x_2)=0$. We aim to reconstruct the unknown pollution source $\lambda_1\delta_{S_1}+\lambda_2\delta_{S_2} $, the unknown boundary data   and the corresponding pollutant concentration field $u$. The ground truths are $\lambda_1=\lambda_2=1$, $S_1=(0.3,0.7)$, $S_2=(0.7,0.3)$ and $g(x_1,x_2)=0.1(\sin(\pi x)-\sin(\pi y))$. The measurement domain is $\omega=(0.4,0.6)^2$.
\end{example}
We denote the coordinates of the two point sources by $\bm{\mu}$ and their intensities by $\bm{\lambda}=(\lambda_1,\lambda_2)$. For the boundary condition $g(x_1,x_2)$, we first parameterize it by arc length as $g(s)$, $s\in[0,4)$, with $s=0$ at the origin $(0,0)$ and increasing counterclockwise, and write $g(s)\approx\sum_{i=1}^{N_i} g^i\,\phi_i(s)$, with piecewise linear hat functions $\{\phi_i\}_{i=1}^{N_g}$ on $[0,4)$ ($N_g=40$ in the experiments). The penalty parameters take the values $\lambda_{\bm{\mu}}^k\equiv\text{5e-4}$ and $\lambda_{g}^k=\text{5e-4}\times0.9^{k}$. 

\begin{figure}
     \centering	
     \begin{tabular}{cccc}
		\includegraphics[width=0.20\textwidth]{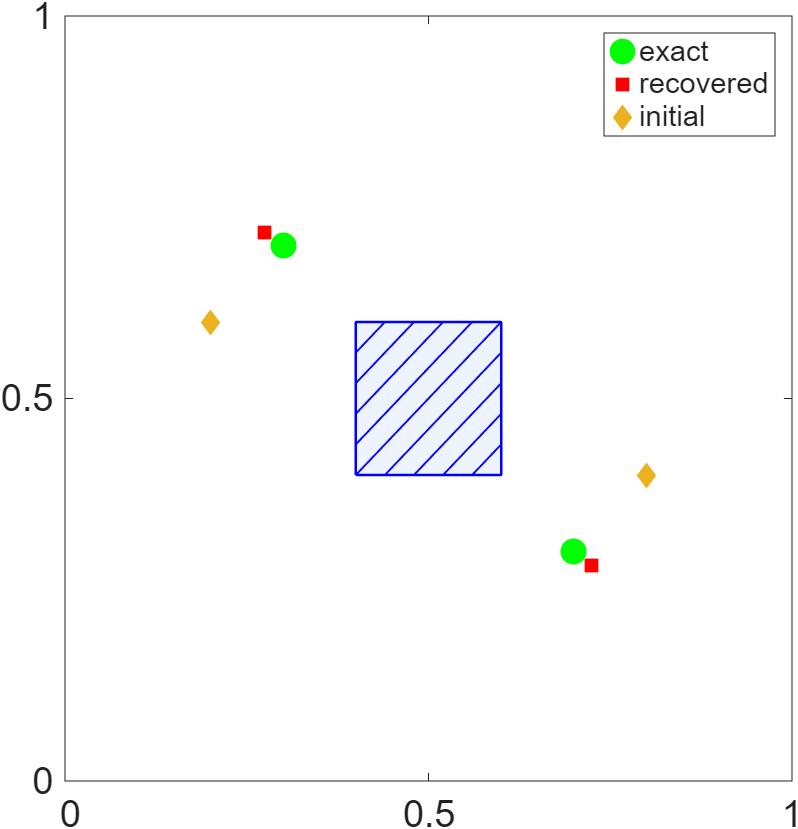}&
		\includegraphics[width=0.20\textwidth]{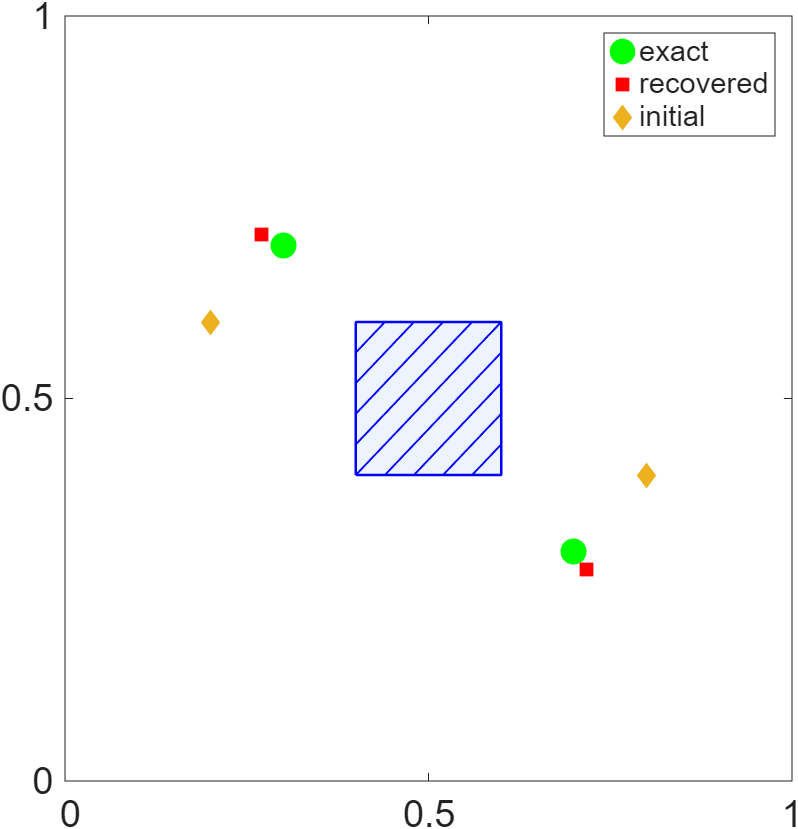}  &
		\includegraphics[width=0.24\textwidth]{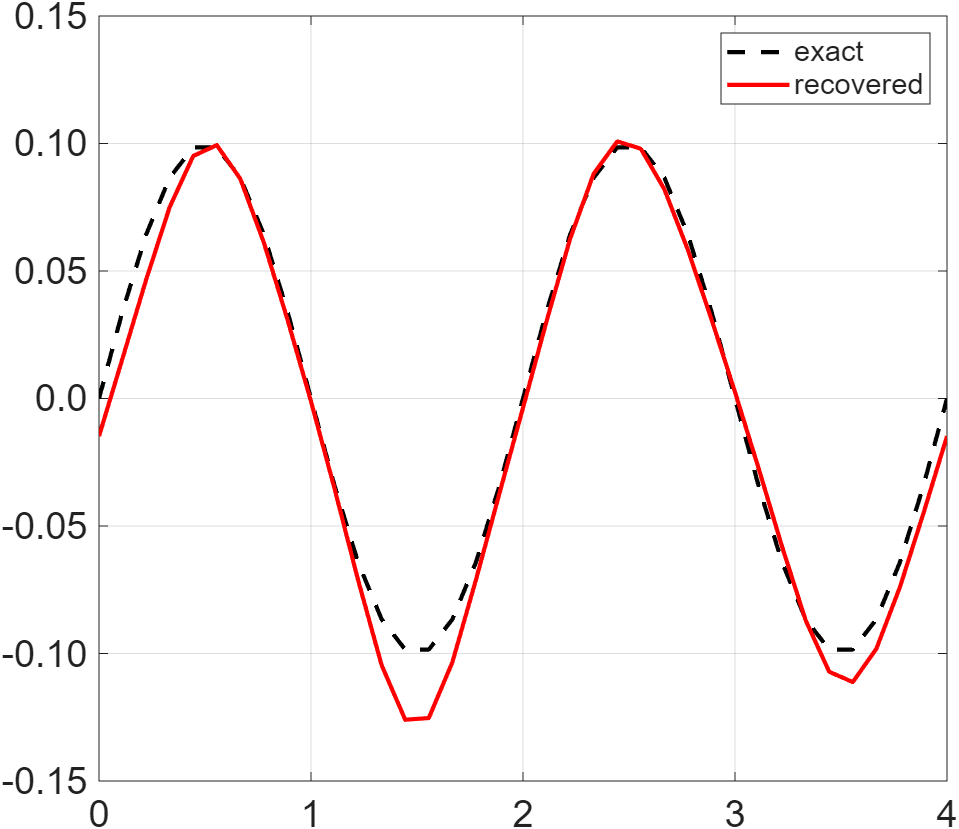} &
        \includegraphics[width=0.24\textwidth]{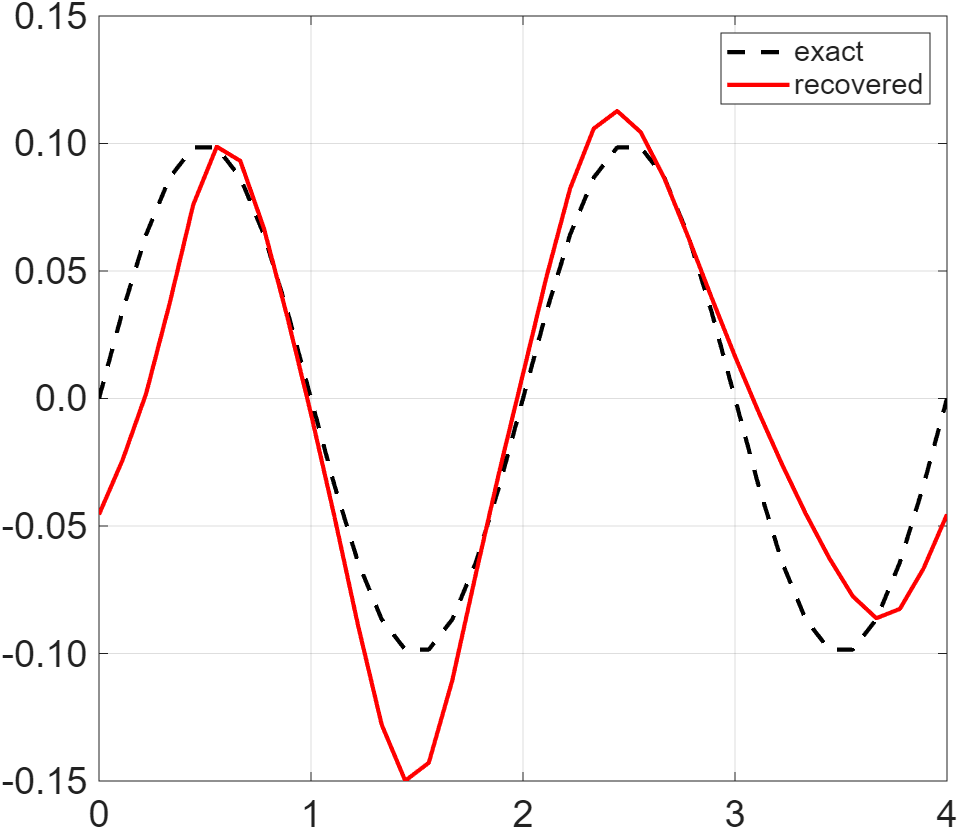}\\ 
        (a) source with $\delta=1\%$  & (b) source with $\delta=5\%$  & (c) $g$ with $\delta=1\%$  & (d)   $g$ with $\delta=5\%$
	\end{tabular}  
    \caption{The reconstructions of pollution sources and boundary data  for Example~\ref{ex:bdry} at different noise levels.  }	\label{Fig:bdry_source}
 \end{figure}

\begin{figure}
     \centering	
     \begin{tabular}{ccc} 
		\includegraphics[width=0.24\textwidth]{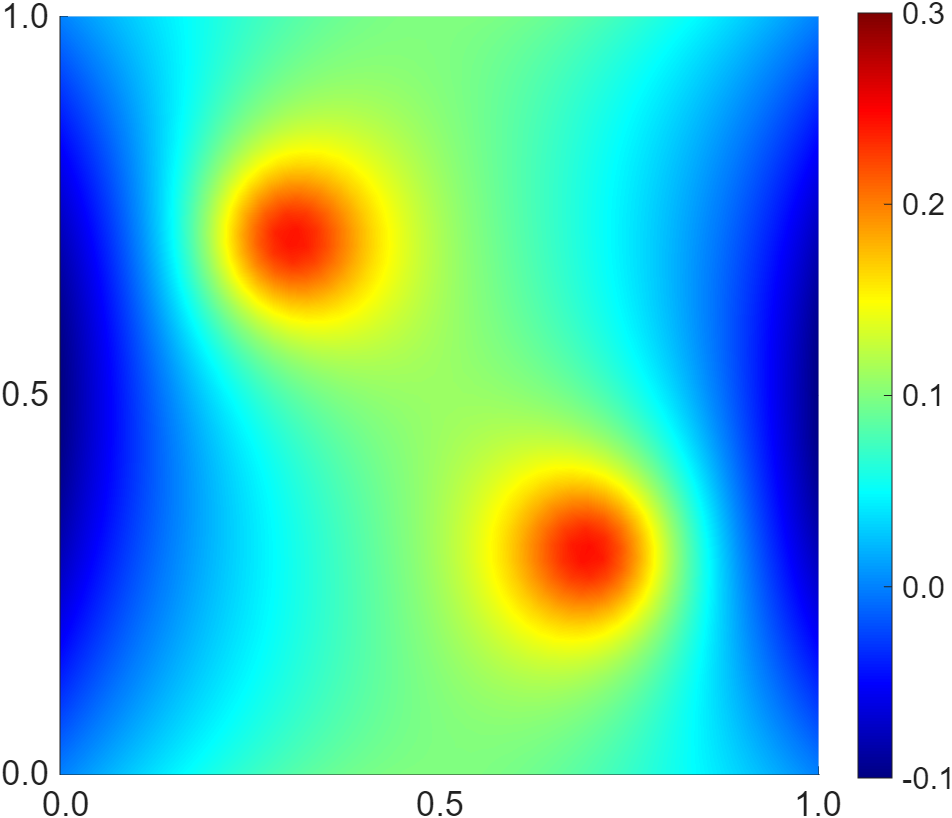}  &
		\includegraphics[width=0.24\textwidth]{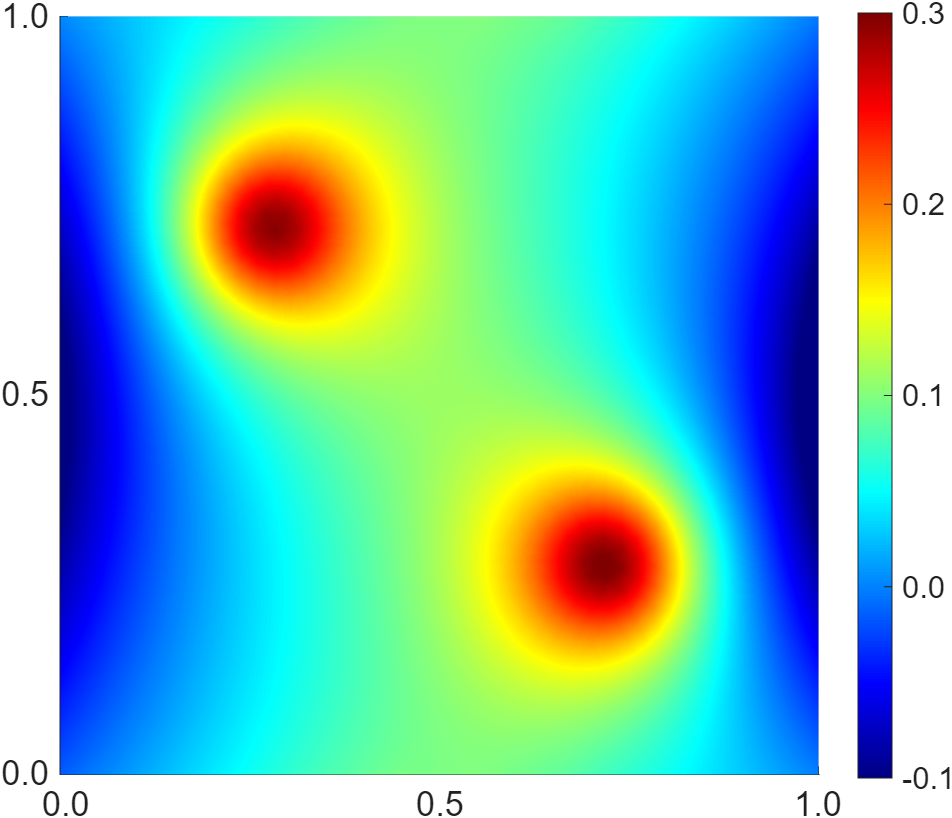} &
        \includegraphics[width=0.24\textwidth]{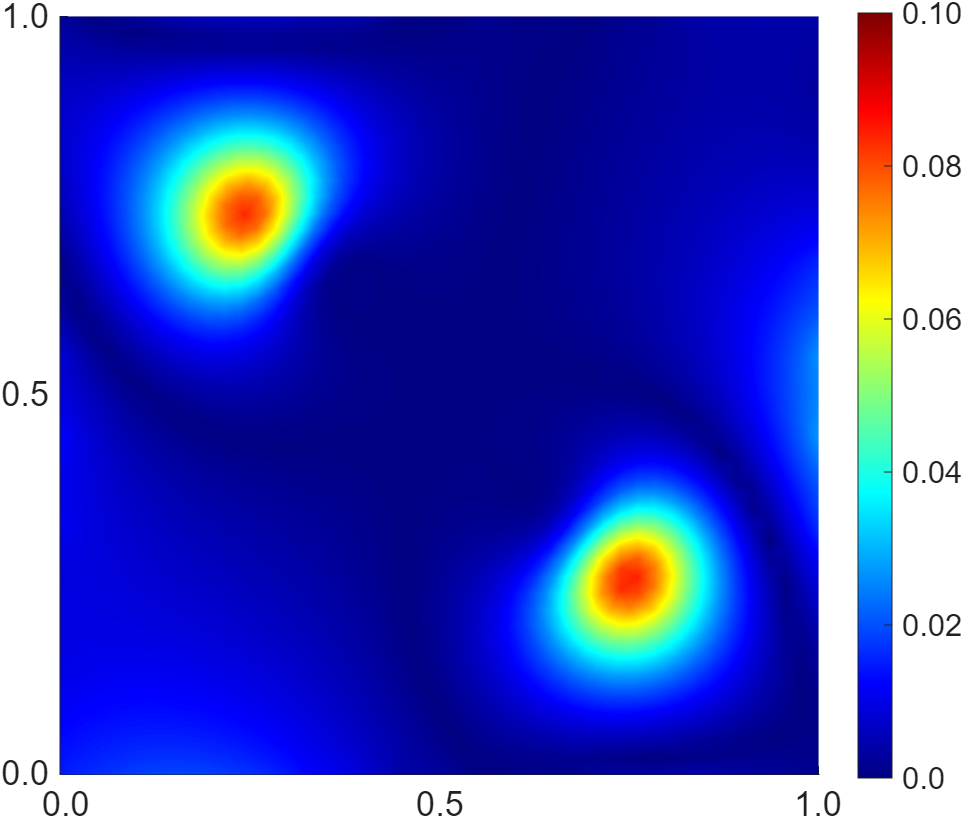} \\
        \includegraphics[width=0.24\textwidth]{pics/example_bdry_ue.png}  &
		\includegraphics[width=0.24\textwidth]{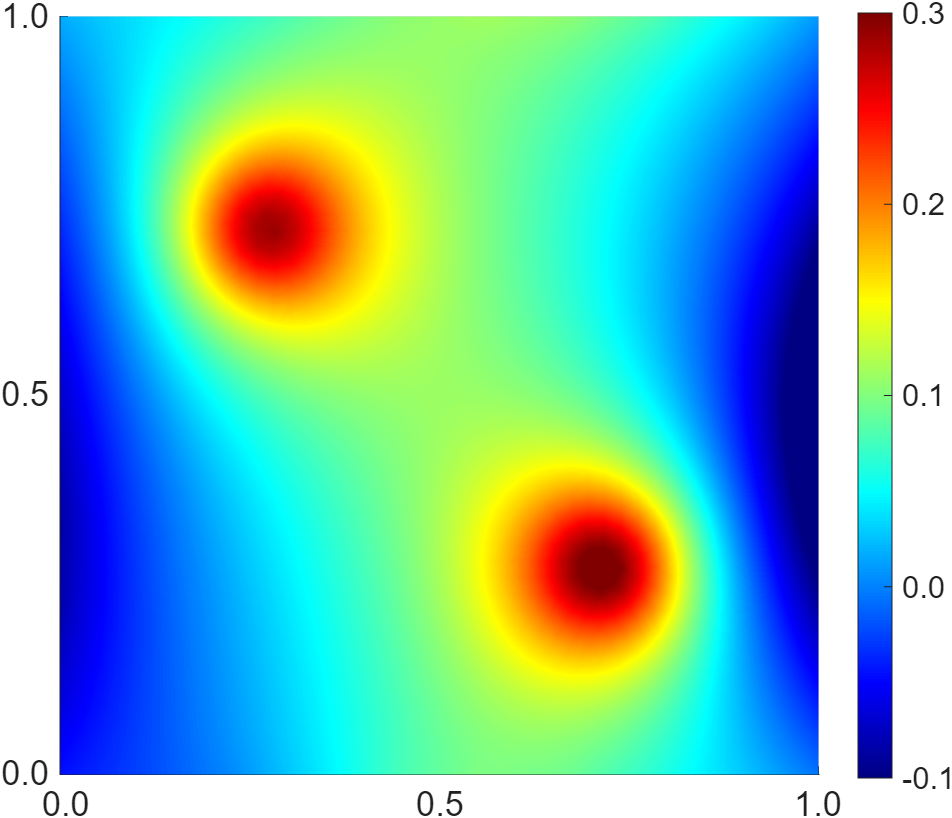} &
        \includegraphics[width=0.24\textwidth]{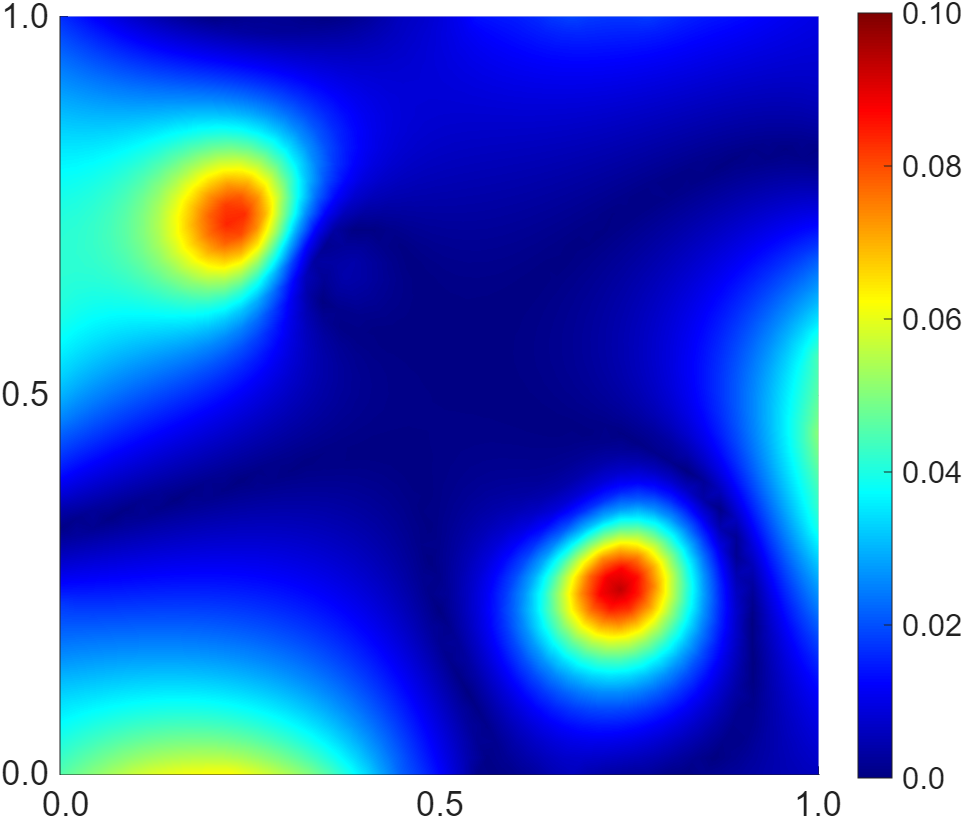} \\
         (a) true $u^\dagger$ & (b) reconstructed $u$ & (c)  error $|u^\dagger-u|$
	\end{tabular}  
    \caption{The reconstructions of pollutant concentration fields $u$ for Example~\ref{ex:bdry} at different noise levels. First row: $1\%$ noise. Second row: $5\%$ noise.}	\label{Fig:bdry_sol}
 \end{figure}

The reconstruction results are shown in Figures~\ref{Fig:bdry_source}, \ref{Fig:bdry_sol} and Table~\ref{tab:error}. We observe that the reconstructions of the pollution sources are relatively accurate and stable. However, the reconstruction of the boundary data deteriorates at high noise levels. This difficulty is inherent to the problem: recovering $g$ from the measurements on $\omega$ amounts to a unique continuation from $\omega$ to the boundary $\partial\Omega$, whose stability deteriorates with the distance.  

\section{Conclusion}

In this article, we consider the problem of determining the global concentration field of a contaminant in groundwater from water extracted by sampling tubes located in the observation region \eqref{omega}. To this end, by exploiting mathematical models describing the diffusion of groundwater pollution, we reformulate this problem in terms of the \textbf{(GUC)} problem. Taking into account both steady-state and unsteady-state diffusion processes, we introduce two different variants, namely \textbf{(SGUC)} and \textbf{(NSGUC)}, of the \textbf{(GUC)} problem. We also establish theoretical results, stated in Theorems \ref{t2}-\ref{t4}, which provide what appears to be the first positive answer to these problems. This theoretical analysis is complemented by numerical simulations illustrating the reconstruction of the concentration field from noisy local measurements, using a regularized Levenberg-Marquardt method. Different numerical experiments illustrate the theory and confirm that the concentration field can be recovered accurately. Our analysis provides a valuable toolbox for  applications to the reconstruction of aquifer hydraulic properties.

This analysis could be extended to investigate finer theoretical properties, such as stability estimates. Alongside such extensions, the precision and performance of the numerical simulations could be further enhanced by incorporating approaches such as rigorous error analysis. We leave these extensions to future investigations.
\section*{Acknowledgment}

The authors would like to warmly thank Val\'erie Mesnage for fruitful discussions on this project and its experimental motivations from problems in hydrology and geoscience, but also for sharing with us the different pictures illustrating this article. The authors are also grateful to the multidisciplinary project FORECAST, conducted jointly by the laboratory of Math\'ematiques Rapha\"el Salem and the Morphodynamique Continentale et C\^oti\`ere laboratory of the University of Rouen, where this work was initiated.

\bibliographystyle{abbrv}
\bibliography{ref.bib}

@article{SS1,
  author = {Saut, J. C. and Scheurer, B.},
  title = {Sur l'unicit{\'e} du probl{\`e}me de Cauchy et le prolongement unique pour des {\'e}quations elliptiques {\`a} coefficients non localement born{\'e}s},
  journal = {J. Diff. Equat.},
  volume = {43},
  year = {1982},
  pages = {28--43}
}

@article{SS,
  author = {Saut, J. C. and Scheurer, B.},
  title = {Unique continuation for some evolution equations},
  journal = {J. Diff. Equat.},
  volume = {66},
  year = {1987},
  pages = {118--139}
}

@article{MS,
  author    = {Mahinthakumar, G. (Kumar) and Sayeed, Mohamed},
  title     = {Hybrid Genetic Algorithm—Local Search Methods for Solving Groundwater Source Identification Inverse Problems},
  journal   = {Journal of Water Resources Planning and Management},
  year      = {2005},
  volume    = {131},
  number    = {1},
  pages     = {45--57},
  publisher = {American Society of Civil Engineers},
  doi       = {10.1061/(ASCE)0733-9496(2005)131:1(45)},
  url       = {https://doi.org/10.1061/(ASCE)0733-9496(2005)131:1(45)}
}

@article{GL,
 author = {Garofalo, Nicola and Lin, Fang-Hua},
 title = {Unique continuation for elliptic operators: {A} geometric-variational approach},
 fjournal = {Communications on Pure and Applied Mathematics},
 journal = {Commun. Pure Appl. Math.},
 issn = {0010-3640},
 volume = {40},
 number = {3},
 pages = {347--366},
 year = {1987},
 language = {English},
 doi = {10.1002/cpa.3160400305},
 zbMATH = {4103553},
 Zbl = {0674.35007}
}

@article{Aliouache2025,
  author = {Aliouache, M. and Jardani, A. and Krimissa, M. and Vu, T. M. and Cabana Lourino, B. and Lorthioy, M.},
  title = {Characterization of the aquifer's hydraulic properties using different pumping test modes: From oscillatory hydraulic tomography to steady state hydraulic tomography},
  journal = {Journal of Hydrology},
  volume = {662},
  pages = {133897},
  year = {2025}
}

@book{Bear1972,
  author = {Bear, J.},
  title = {Dynamics of Fluids in Porous Media},
  publisher = {Elsevier},
  address = {New York},
  year = {1972}
}

@article{Yeh2000,
  author = {Yeh, T.-C. J. and Liu, S.},
  title = {Hydraulic tomography: Development of a new aquifer test method},
  journal = {Water Resources Research},
  volume = {36},
  number = {8},
  pages = {2095--2105},
  year = {2000}
}

@article {ElBadiaHaDuong2002pollution_detection_problem,
    AUTHOR = {El Badia, A. and Ha-Duong, T.},
     TITLE = {On an inverse source problem for the heat equation.
              {A}pplication to a pollution detection problem},
   JOURNAL = {J. Inverse Ill-Posed Probl.},
  FJOURNAL = {Journal of Inverse and Ill-Posed Problems},
    VOLUME = {10},
      YEAR = {2002},
    NUMBER = {6},
     PAGES = {585--599},
      ISSN = {0928-0219,1569-3945},
   MRCLASS = {35R30 (35K05 35K20 80A23)},
  MRNUMBER = {1967440},
MRREVIEWER = {Peter\ G.\ Danilaev},
       DOI = {10.1515/jiip.2002.10.6.585},
       URL = {https://doi-org.ezproxy.lb.polyu.edu.hk/10.1515/jiip.2002.10.6.585},
}

@article{gong2026identification,
  TITLE = {{Identification of a Point Source in the Heat Equation from Sparse Boundary Measurements }},
  AUTHOR = {Gong, Fangyu and Jin, Bangti and Kian, Yavar and Liu, Sizhe},
  URL = {https://hal.science/hal-05545664},
   JOURNAL = {SIAM J. Math. Anal.},
  FJOURNAL = {SIAM Journal on Mathematical Analysis},
  PUBLISHER = {{Society for Industrial and Applied Mathematics}},
  YEAR = {2026 in press},
  HAL_ID = {hal-05545664},
  HAL_VERSION = {v1},
}

@article{Bu1,
 author = {Burman, Erik and Delay, Guillaume and Ern, Alexandre},
 title = {The unique continuation problem for the heat equation discretized with a high-order space-time nonconforming method},
 fjournal = {SIAM Journal on Numerical Analysis},
 journal = {SIAM J. Numer. Anal.},
 issn = {0036-1429},
 volume = {61},
 number = {5},
 pages = {2534--2557},
 year = {2023},
 language = {English},
 doi = {10.1137/22M1508637},
 zbMATH = {7764838},
 Zbl = {1549.65349}
}

@article{Nag,
  title={Compactifications of unstable N{\"o}beling spaces},
  author={Nag{\'o}rko, Artur},
  journal={arXiv preprint arXiv:1712.03181},
  year={2017}
}

@book{LM1,
 author = {Lions, J. L. and Magenes, E.},
 title = {Non-homogeneous boundary value problems and applications. {Vol}. {II}. {Translated} from the {French} by {P}. {Kenneth}},
 fseries = {Grundlehren der Mathematischen Wissenschaften},
 series = {Grundlehren Math. Wiss.},
 issn = {0072-7830},
 volume = {182},
 year = {1972},
 publisher = {Springer, Cham},
 language = {English},
 zbMATH = {3360568},
 Zbl = {0227.35001}
}

@book{GT,
 author = {Gilbarg, David and Trudinger, Neil S.},
 title = {Elliptic partial differential equations of second order. 2nd ed},
 fseries = {Grundlehren der Mathematischen Wissenschaften},
 series = {Grundlehren Math. Wiss.},
 issn = {0072-7830},
 volume = {224},
 year = {1983},
 publisher = {Springer, Cham},
 language = {English},
 zbMATH = {3895476},
 Zbl = {0562.35001}
}

@book{Hir,
  title={Differential Topology},
  author={Hirsch, Morris W.},
  series={Graduate Texts in Mathematics},
  volume={33},
  year={1976},
  publisher={Springer-Verlag},
  address={New York-Berlin},
  isbn={978-0-387-90148-0}
}

@Article{w12092415,
AUTHOR = {Jamshidi, Azade and Samani, Jamal Mohammad Vali and Samani, Hossein Mohammad Vali and Zanini, Andrea and Tanda, Maria Giovanna and Mazaheri, Mehdi},
TITLE = {Solving Inverse Problems of Unknown Contaminant Source in Groundwater-River Integrated Systems Using a Surrogate Transport Model Based Optimization},
JOURNAL = {Water},
VOLUME = {12},
YEAR = {2020},
NUMBER = {9},
ARTICLE-NUMBER = {2415},
URL = {https://www.mdpi.com/2073-4441/12/9/2415},
ISSN = {2073-4441},
DOI = {10.3390/w12092415}
}

@article{theis,
  title={The relation between the lowering of the piezometric surface and the rate and duration of discharge of a well using ground-water storage},
  author={Theis, Charles V.},
  journal={Eos, Transactions American Geophysical Union},
  volume={16},
  number={2},
  pages={519--524},
  year={1935},
  publisher={Wiley Online Library}
}

@book{huyak,
  title={Computational Methods in Subsurface Flow},
  author={Huyakorn, Peter S. and Pinder, George F.},
  year={1983},
  publisher={Academic Press},
  address={New York}
}

@book{sch,
  title={The Physics of Flow Through Porous Media},
  author={Scheidegger, Adrian E.},
  edition={3rd},
  year={1974},
  publisher={University of Toronto Press},
  address={Toronto}
}

@article{Bu2,
 author = {Burman, Erik and Oksanen, Lauri},
 title = {Data assimilation for the heat equation using stabilized finite element methods},
 fjournal = {Numerische Mathematik},
 journal = {Numer. Math.},
 issn = {0029-599X},
 volume = {139},
 number = {3},
 pages = {505--528},
 year = {2018},
 language = {English},
 doi = {10.1007/s00211-018-0949-3},
 zbMATH = {6904222},
 Zbl = {1412.65111}
}

@misc{DCK,
 author = {Ulysse Dalmasso and Siyu Cen and Yavar Kian},
 title = {Recovery of a {Measure}-valued {Source} in the {Heat} {Equation} from {Sparse} {Boundary} {Measurements}},
 year = {2026},
 howpublished = {Preprint, {arXiv}:2607.15841 [math.{AP}] (2026)},
 url = {https://arxiv.org/abs/2607.15841},
 arXiv = {arXiv:2607.15841}
}

@article{GZZ,
 author = {Gu, Qiling and Zhang, Wenlong and Zhang, Zhidong},
 title = {Determine the point source of the heat equation with sparse boundary measurements},
 fjournal = {SIAM Journal on Applied Mathematics},
 journal = {SIAM J. Appl. Math.},
 issn = {0036-1399},
 volume = {85},
 number = {5},
 pages = {2337--2354},
 year = {2025},
 language = {English},
 doi = {10.1137/25M1725620},
 zbMATH = {8113966},
 Zbl = {1577.35435}
}

@article {Andrle2012moving_pollution_sources,
    AUTHOR = {Andrle, M. and El Badia, A.},
     TITLE = {Identification of multiple moving pollution sources in surface
              waters or atmospheric media with boundary observations},
   JOURNAL = {Inverse Problems},
  FJOURNAL = {Inverse Problems. An International Journal on the Theory and
              Practice of Inverse Problems, Inverse Methods and Computerized
              Inversion of Data},
    VOLUME = {28},
      YEAR = {2012},
    NUMBER = {7},
     PAGES = {075009, 22 pp.},
      ISSN = {0266-5611,1361-6420},
   MRCLASS = {86A10 (35K20 35R30 86A05)},
  MRNUMBER = {2946797},
MRREVIEWER = {Sven\ Ivansson},
       DOI = {10.1088/0266-5611/28/7/075009},
       URL = {https://doi-org.ezproxy.lb.polyu.edu.hk/10.1088/0266-5611/28/7/075009},
}

@article{GK, 
   author = {Gurarslan, Gurhan and Karahan, Halil},  
    title = {Solving inverse problems of groundwater-pollution-source identification using a differential evolution algorithm}, 
   journal = {Hydrogeol. J.}, 
    volume = {23},
      year = {2015}, 
     pages = {1109--1119},
       doi = {https://doi.org/10.1007/s10040-015-1256-z},
}

@article{Zhu2005,
  author = {Zhu, J. and Yeh, T.-C. J.},
  title = {Characterization of aquifer heterogeneity in transect using transient hydraulic tomography},
  journal = {Water Resources Research},
  volume = {41},
  number = {5},
  pages = {W05008},
  year = {2005}
}

@book {Thomee:2006,
    AUTHOR = {Thom\'{e}e, Vidar},
     TITLE = {Galerkin {F}inite {E}lement {M}ethods for {P}arabolic {P}roblems},
   EDITION = {Second},
 PUBLISHER = {Springer-Verlag, Berlin},
      YEAR = {2006},
     PAGES = {xii+370},
      ISBN = {978-3-540-33121-6; 3-540-33121-2},
   MRCLASS = {65-02 (65M15 65M60)},
  MRNUMBER = {2249024},
}

@article {Marquardt:1963,
    AUTHOR = {Marquardt, Donald W.},
     TITLE = {An algorithm for least-squares estimation of nonlinear
              parameters},
   JOURNAL = {J. Soc. Indust. Appl. Math.},
  FJOURNAL = {Journal of the Society for Industrial and Applied Mathematics},
    VOLUME = {11},
      YEAR = {1963},
     PAGES = {431--441},
      ISSN = {0368-4245},
   MRCLASS = {62.20},
  MRNUMBER = {153071},
MRREVIEWER = {M.\ Atiqullah},
}

@article {Levenberg:1944,
    AUTHOR = {Levenberg, Kenneth},
     TITLE = {A method for the solution of certain non-linear problems in
              least squares},
   JOURNAL = {Quart. Appl. Math.},
  FJOURNAL = {Quarterly of Applied Mathematics},
    VOLUME = {2},
      YEAR = {1944},
     PAGES = {164--168},
      ISSN = {0033-569X,1552-4485},
   MRCLASS = {65.0X},
  MRNUMBER = {10666},
MRREVIEWER = {T.\ E.\ Sterne},
       DOI = {10.1090/qam/10666},
       URL = {https://doi.org/10.1090/qam/10666},
}

@article{Hecht:2012,
  AUTHOR = {Hecht, F.},
  TITLE = {New development in FreeFem++},
  JOURNAL = {J. Numer. Math.},
  FJOURNAL = {Journal of Numerical Mathematics},
  VOLUME = {20}, YEAR = {2012},
  NUMBER = {3-4}, PAGES = {251--265},
  ISSN = {1570-2820},
  MRCLASS = {65Y15},
  MRNUMBER = {3043640},
  URL = {https://freefem.org/}
}

@article{Hughes:1982,
 author = {Brooks, Alexander N. and Hughes, Thomas J. R.},
 title = {Streamline upwind/{Petrov}-{Galerkin} formulations for convection dominated flows with particular emphasis on the incompressible {Navier}-{Stokes} equations},
 fjournal = {Computer Methods in Applied Mechanics and Engineering},
 journal = {Comput. Methods Appl. Mech. Eng.},
 issn = {0045-7825},
 volume = {32},
 pages = {199--259},
 year = {1982} 
}
 \end{document}